%% file: saw_etds.tex
\documentclass{article}

\usepackage{graphicx}
\usepackage{amsmath,amssymb, amsthm}
\usepackage{subfig}
\usepackage{float}
\usepackage{standalone}
\usepackage{tikz}
\usepackage{geometry}
\usepackage{cite}
\usetikzlibrary{calc, intersections, decorations.pathreplacing, positioning}
\newtheorem{theorem}{Theorem}[section]

\newtheorem{lemma}[theorem]{Lemma}
\newtheorem{proposition}{Proposition}

\theoremstyle{definition}

\newtheorem{remark}{Remark}

\title{\textbf{Chaos in Saw Map}}
\author{Nikita Begun$^{* \dagger \ddagger}$, Pavel Kravetc$^{\S}$ and Dmitrii Rachinskii$^{\S}$\\ \\
	\small $^*$ Institut f\"ur Mathematik, Freie Universit\"at Berlin, Arnimallee 3 D - 14195, \\
	\small Berlin, Germany \\
	\small $^\dagger$  Saint Petersburg State University, University Embankment, 7/9, \\
	\small Saint-Petersburg, Russia, 199034\\
	\small $^\ddagger$ People's Friendship University of Russia (RUDN University), Miklouho-Maclay St, 6 \\ 				\small Moscow, Russia, 117198\\
	\small (email: nikitabegun88@gmail.com)\\
	\small $^\S$ University of Texas at Dallas, 800 W Cambell Rd, Richardson, TX, USA, 75080\\
	\small (email: pavel.kravetc@utdallas.edu) \\
	\small (email: dmitry.rachinskiy@utdallas.edu)
}

\begin{document}

%% Enter the first author's name and address:
%\centerline{\scshape Nikita Begun$^*$}
%\medskip
%{\footnotesize
%% please put the address of the first author
% \centerline{Institut f\"ur Mathematik}
%   \centerline{Freie Universit\"at Berlin}
%   \centerline{Arnimallee 3 D - 14195, Berlin}
%} % Do not forget to end the {\footnotesize by the sign }
%
%\medskip
%
%\centerline{\scshape Pavel Kravetc and Dmitrii Rachinskii}
%\medskip
%{\footnotesize
% % please put the address of the second  and third author
% \centerline{Department of Mathematical Sciences}
%   \centerline{The University of Texas at Dallas}
%   \centerline{Richardson, TX 75080, USA}
%}
%
%

\maketitle

\begin{abstract}
	We consider dynamics of a scalar piecewise linear  "saw map" with infinitely many linear segments.
	In particular, such maps are generated as a Poincar\'e map of simple two-dimensional discrete time piecewise linear systems
	involving a saturation function. Alternatively, these systems can be viewed as
	a feedback loop with the so-called stop hysteresis operator. We analyze chaotic sets and attractors of
	the "saw map" depending on its parameters.
\end{abstract}

\section{Introduction}

Piecewise linear (PWL) and piecewise smooth (PWS) systems serve as a modeling framework for applications involving friction, collision, sliding, intermittently constrained systems and processes with switching components. Examples include stick-slip motion in mechanical systems, impact oscillators, switching electronic circuits (such as DC/DC power converters and relay controllers), hybrid dynamics in control systems as well as models of economics and finance \cite{2003, Champneys}. Low-dimensional nonsmooth maps can also appear as Poincar\'{e} return maps of smooth flows which show chaotic dynamics \cite{Guckenheimer,Zaks}. 

Since the second half of the 1990s, multiple analytic tools have been developed in order to address distinctive scenarios, which are not observed in smooth dynamical systems but are unique to piecewise smooth systems \cite{Survey,Simpson,Kunze}\footnote{Some of these methods are based on important results that had already been obtained before \cite{Andronov, Feigin, Leonov, Filippov}.}. For example, these scenarios include robust chaos \cite{rc}; various types of discontinuity-induced bifurcations (related to sliding, chattering, grazing and corner collision phenomena) which occur when an invariant set collides with a switching surface; and, border-collision bifurcations of the associated Poincar\'{e} maps with their normal forms \cite{NY, meist,meiss,gardini,sushko}.

One-dimensional maps serve as important prototype models, which help understand dynamics of higher-dimensional systems. Further, one-dimensional dynamics has developed into a subject in its own right \cite{demelo}. In particular, bifurcation scenarios and symbolic dynamics specific to PWS and PWL maps have been intensively explored in the last two decades \cite{Gallery}. However, they are still understood to a lesser extent than dynamics of smooth maps, and the PWS theory is far from being complete. Even the skew tent map, which is a simple variation of the classical tent map, produces a rich variety of dynamical scenarios which have not been completely described yet (see survey \cite{SuAvr} for the state-of-the-art results).

In this paper, we consider dynamics of a one-dimensional PWL map, which has infinitely many local maximum and minimum points accumulating near an essential discontinuity point, see Fig.~\ref{fig1}. We call it a ``saw map". As a matter of fact, this map represents a reduction of a two-dimensional PWL map including a linear term and the simple saturation PWL function (see \eqref{saturation} below) to a one-dimensional Poincar\'e map \cite{siamds}. 
The objective of this work is to analyze chaotic repelling and attracting sets, including robust chaos, for a general class of such
one-dimensional maps.
%We establish possible structures of attractors, including robust chaos, which depend on certain geometric parameters of the saw map. 
The ``saw map" can have multiple attractors embedded into a number of invariant intervals. One of such intervals contains infinitely many local maximum and minimum points of the map, while the restriction of the map to any other invariant interval is a skew tent map.  
The paper is organized as follows. In the main Section 2, we consider dynamics of the ``saw map''
depending on its parameters.
In particular, we obtain a characterization 
% description
of the attractor of the skew tent map, which to the best of our knowledge is new (see Remark \ref{remj1} to formula \eqref{Lambda} of Theorem \ref{Jan1}). In Section 3, we discuss the implications of our results for the 2-dimensional PWL map with the saturation nonlinearity. This map has been related to a set of macroeconomic models with sticky inflation and to discrete time systems with dry friction in \cite{siamds}.
\section{Main results}\label{main}

There are several definitions of a chaotic invariant set. We use the definition of R. Devaney \cite{devaney}.

A closed invariant set $A$ for the map $f$ is chaotic if:
\begin{itemize}

\item (density of periodic orbits) periodic points of $f$ are dense in $A$;

\item (sensitivity to initial conditions) there is a $\beta>0$ such that for any $x\in A$ and any $\epsilon>0$ there is a $y\in A$ with $|x-y|<\epsilon$ and a $k$ such that $|f^{k}(x)-f^{k}(y)|>\beta$;

\item (topological mixing or transitivity) %for any pair of non-empty open sets $B,C \subset A$ there exists an $n$ such that $f^{n}(B)\cap C\neq \emptyset$.
$A$ contains a dense orbit of $f$.

\end{itemize}

In what follows, we use the notation $A\subset B$ if $A\subseteq B$, $A\ne B$. By $\overline A$ we denote the closure of a set $A$.

By $|I|$ we denote the length of an interval $I$.
%Sometimes, instead of topological mixing, a stricter condition is required:
%\begin{itemize}

%\item $A$ contains a dense orbit of $f$.

%\end{itemize}

%This condition implies topological mixing.

\begin{figure}[h]
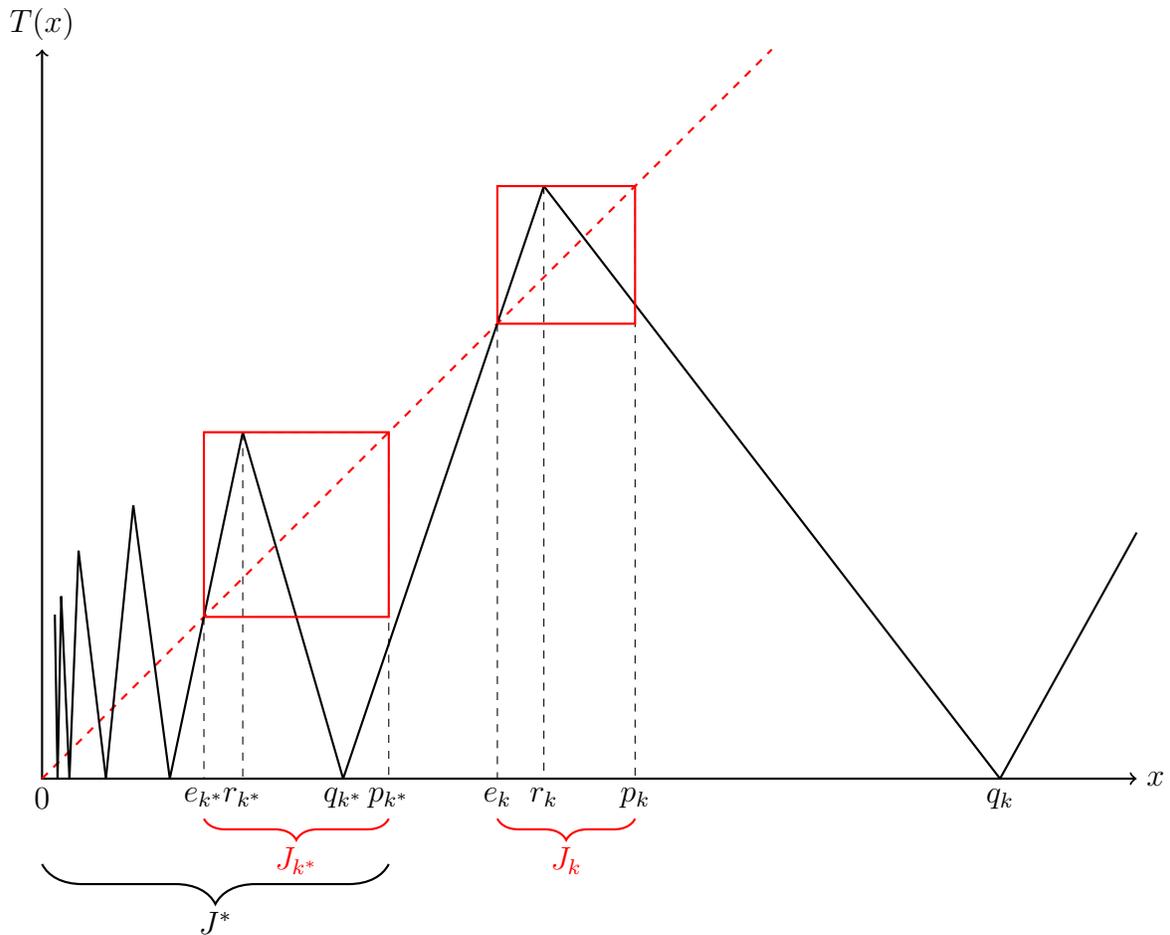

\begin{center}
\includestandalone[width=\textwidth]{saw_pic}
\caption{Graph of the map $T$.} \label{fig1}
\end{center}
\end{figure}

Let sequences $p_k,q_k,r_k$ satisfy
\begin{equation}\label{pqr}
\begin{array}{ll}
r_0>q_1>r_1>q_2>r_2>\cdots>0,&\quad q_k,r_k\to 0,\\
p_0>p_1>p_2>p_3>\cdots>0,& \quad p_0<r_0, \quad p_k>r_k \ \ {\rm for\ all}\ \ k\geq1.
%,&\qquad p_k\to p_*.
\end{array}
\end{equation}
Denote $J=[0,r_0]$ and consider a ``saw map" $T:J\to \mathbb{R}$ defined by the following properties:
\begin{itemize}
\item $T(0)=T(q_k)=0$ for all 	$k\ge 1$;

\item $T(r_k)=p_k$ for all 	$k\ge 0$;

\item $T$ is linear on each of the intervals $[q_{k+1},r_k]$, $[r_k,q_k]$, $k\ge1$, and $[q_1,r_0]$,
\end{itemize}
see Fig.~\ref{fig1}.
These properties and the fact that $p_0<r_0$ imply that $T(J) \subset J$. Also these properties imply that $T$ has a unique fixed point $e_k$ in each interval $(q_{k+1},r_k)$, $k \geq 1$
and a unique fixed point $\hat e_k$ in each interval $(r_k,q_k)$, $k\geq 1$.

Note that $T$ is piecewise linear and continuous on every segment $[a,b]\subset J\setminus\{0\}$.
If $p_*:=\lim_{k\to\infty} p_k=0$, then $T$ is continuous on $J$.
On the other hand, if $p_*>0$, then $T$ has a discontinuity at zero.
% Consider $k_0$ such that $T(r_{k_0})\geq r_{k_0}$, but $T(r_{k_0-1})<r_{k_0-1}$, see Fig. \ref{fig1}. Denote by $J$ the invariant segment $[0,T(r_{k_0})]$.
Denote by
\[
\alpha_k:=p_k/(r_k-q_{k+1}),\qquad \beta_k:=p_k/(q_k-r_k)
\]
the absolute value of the slope of the graph of $T$ on the intervals $[q_{k+1},r_k]$ and $[r_k,q_k]$, respectively.
The assumption $p_k>r_k$, $k\geq 1$ implies that
\[
\alpha_k>1 \quad
{\rm for\ all} \quad
k\geq 1.
\]

 We assume additionally that
\begin{itemize}
\item $\alpha_k$ and $\beta_k$ are increasing sequences;

\item
%There exists a $\bar{k}$ such that $T(r_{\bar{k}})\leq e_{\bar{k}-1}$. Then,
there exists  $k^*\ge 1$ such that $p_k>e_{k-1}$ for $k\ge k^*+1$ and if $k_*\ge 2$ then
$p_k < e_{k-1}$ for $2 \leq k\le k^*$;

\item for $1\leq k\leq k^*-1$,
\begin{equation}\label{star}
q_k+\frac{e_k}{\alpha_{k-1}}>p_k.
\end{equation}

\end{itemize}

The last condition is technical. It allows us to simplify the results by shortening the list of possible dynamical scenarios, which are, however, similar to each other.

%Under the assumptions of Lemma \ref{l2},
From the definition of $k^*$ it follows that the segment
\[
J^*=[0,p_{k^*}]\subset J
\]
is invariant for $T$. A proof of this fact is given in the next section. %If $\bar{k}$ such that $T(r_{\bar{k}})\leq e_{\bar{k}-1}$ does not exists then we set $k^*=k_0,$ $J^*=J$.
%If $k^*\ge 2$,
Consider also the segments
\begin{equation}\label{JG}
J_k=[e_k,p_k], %\quad k\le k^*;
\qquad
G_k=[T^2(r_k),T(r_k)]=[T(p_k),p_k],\quad 1 \leq k\le k^*,%k<k^*.
\end{equation}
see Fig.~\ref{fig1}.
%Obviously, $J^*$ and $J_k$ belong to $J$.
%For $k=k_0,\ldots,k^*$, let us d
%Note that the sequences $\alpha_k$, $\beta_k$ increase with $k$.
%Recall that $\alpha_k>1$, $\alpha_k>\beta_k$ for all $k$.
These segments are well-defined. Indeed, $e_k<r_k<p_k$, $k\geq 1$ by definition of $T$. Also, since $p_k<e_{k-1}$ for $2 \leq k\le k^*$ and $r_0>p_0$, for $1 \leq k\le k^*$ we have $p_k>T(p_k)$.

\begin{remark}\label{remslope2}
Denote %there exist $f_k$ such that
\begin{equation}\label{fk}
f_k=\min\{x> e_k:\ T(x)=e_k\},\qquad 1 \leq k\le k^*.
\end{equation}
	It is easy to see that if $\alpha_k^{-1}+\beta_k^{-1}\geq 1$, then $f_k\geq p_k$ and the segment $J_k$ is invariant,
	while if $\alpha_k^{-1}+\beta_k^{-1}<1$, then $f_k<p_k$ and the segment $J_k$ is not invariant.
	%\end{remark}
	%\begin{remark}\label{remslope2}
	%It is also easy to see that $\alpha_k^{-1}+\beta_k^{-1}<1$ if $J_k$ is not invariant while $\alpha_k^{-1}+\beta_k^{-1}\geq 1$ if $J_k$ is invariant. In particular, since $T(r_k)>e_{k-1}>f_k$ for $k\geq k^*+1$ then $%\alpha_k^{-1}+\beta_k^{-1}<1$ for $k\geq k^*+1$.
\end{remark}

\begin{theorem}\label{T2}
{\em  If $\alpha_{k^*}^{-1}+\beta_{k^*}^{-1}<1$, then the segment $J^*$ is a chaotic invariant set.
  If $\alpha_{k^*}^{-1}+\beta_{k^*}^{-1}\ge 1$, then the segment $J^*$ contains a chaotic invariant Cantor set $\Sigma$. Furthermore, all trajectories from $J^*\backslash\Sigma$ reach
  the invariant segment $J_{k^*}\subset J^*$
  after a finite number of iterations.}
\end{theorem}

In the case $\alpha_{k^*}^{-1}+\beta_{k^*}^{-1}\geq1$, dynamics on the invariant interval $J_{k^*}$ is described in the next theorem.

\begin{theorem}
	\label{Jan1}
	{\em
		If $\beta_k<1$ for some $k\le k^*$ (see Fig.~\ref{fig:types}(a)), then the interval $J_k$ is invariant and $J_k\setminus \{e_k\}$ belongs to the basin of attraction of the stable fixed point $\hat e_k\in J_k$.

		If $\beta_k>1$ and $\alpha_{k}^{-1}+\beta_{k}^{-1}> 1$ for some $k\le k^*$ (see Fig.~\ref{fig:types}(b)), then the interval $J_k$ is also invariant. Further, the interval $G_k\subset J_k$ defined in \eqref{JG} is invariant and contains a chaotic invariant set
		\begin{equation}\label{Lambda}
		\Lambda_k=\bigcup_{i=1}^{2^N} \overline{A_i},
		\end{equation}
		 where the $A_i$ are non-intersecting open intervals.
		 The complement $G_k\setminus \Lambda_k$ contains $N$ unstable periodic orbits $O_i$ of periods $1,\ldots,2^{N-1}$. Furthermore, for any
		 $x\in G_k\setminus (\Lambda_k\cup O_1\cup\cdots\cup O_{N})$
		 there is an $n_0$ such that $T^n(x)\in \Lambda_k$ for $n\ge n_0$,
		 and for any $y\in J_k\setminus (G_k\cup\{e_k\})$ there is an $n_1$ such that
		 $T^n(y)\in G_k$ for $n\ge n_1$.
		
		 Finally, if $\alpha_{k}^{-1}+\beta_{k}^{-1}<1$ for some $1\le k\le k^*-1$ (see Fig.~\ref{fig:types}(c)), then the interval $J_k$ contains a chaotic invariant Cantor set $\Lambda_k$,
	%with a dense orbit,
	dynamics on this set is conjugate to the left shift, and for every $x\in J_k\setminus \Lambda_k$
	there is an $n_2$ such that $T^n(x)\in [0,e_k)$ for $n\ge n_2$.
		 }
		 
 	\begin{figure}[h]
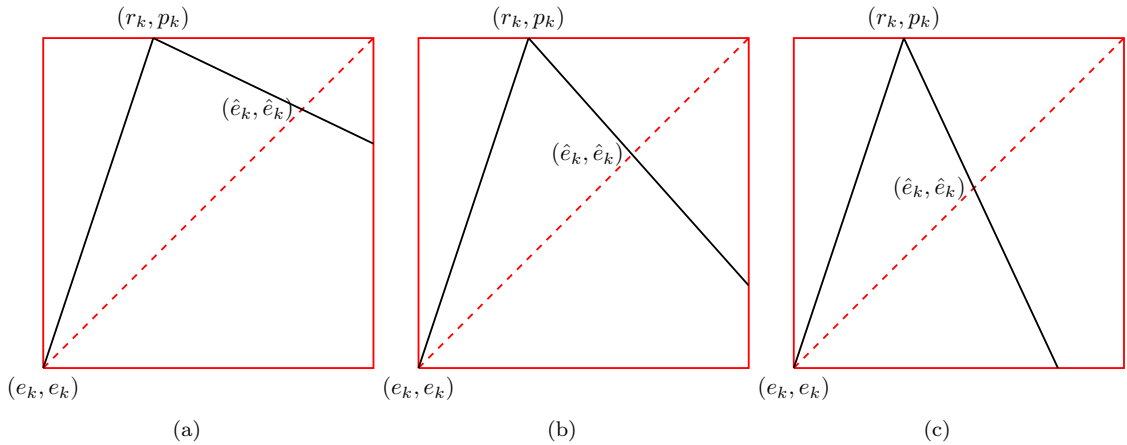

 	\centering
 		\subfloat[\label{fig:type1}]{\includestandalone[width=0.32\textwidth]{type1}}
 		\subfloat[\label{fig:type2}]{\includestandalone[width=0.32\textwidth]{type2}}
 		\subfloat[\label{fig:type3}]{\includestandalone[width=0.32\textwidth]{type3}}
 	\caption{Possible shapes of the graph of the map $T$ on an interval $J_k$. (a) $\beta_k < 1$; (b) $\beta_k > 1$ and $\alpha^{-1}_k + \beta^{-1}_k \ge 1$; (c)  $\beta_k > 1$ and $\alpha^{-1}_k + \beta^{-1}_k < 1$.\label{fig:types}}
 	\end{figure}
\end{theorem}

\begin{remark}\label{remj1}
%	As shown in the proof below, there is a numbering of the intervals $A_i$ such that
%	\[
%	T(A_i)=A_{i+1},\quad i=1,\ldots 2^N-1;\qquad T(A_{2^N})=A_1.
%	\]
As shown in the proof below, 
in the case $\beta_k>1$, $\alpha_{k}^{-1}+\beta_{k}^{-1}> 1$ (see Fig.~\ref{fig:types}(b)),
there is a permutation $\sigma_{2^N}=\bigl( \sigma^1_{2^N}, \sigma^2_{2^N}, ... , \sigma^{2^N}_{2^N}\bigr)$ of the natural numbers from $1$ to $2^N$ and the numbering of the intervals $A_i\subset J_k$ such that 
\[
T(A_i)=A_{i+1},\quad i=1,\ldots 2^N-1;\qquad T(A_{2^N})=A_1,
\]
and the interval $A_{\sigma^k_{2^N}}$ lies to the left of interval $A_{\sigma^{k+1}_{2^N}}$, $k=1,...,2^N-1$. Further, the permutation 
$\sigma_{2^{N}}$ %is obtained from permutation $\sigma_{2^N}$ by the following procedure
is defined by the following inductive (in $N$) formulas:
\begin{equation}
\label{sigma}
\sigma^{2^N+i}_{2^{N+1}}=2\sigma^{i}_{2^N}-1, \qquad
\sigma^{i}_{2^{N+1}}=\sigma^{2^{N+1}+1-i}_{2^{N+1}}+1, \qquad i=1,...,2^N.
\end{equation}		
		\end{remark}
		
		\begin{remark}\label{remj2}
		If $\alpha_{k}^{-1}+\beta_{k}^{-1}=1$, then $T$ is a tent map on $J_k$, and the whole segment $J_k$
		is a chaotic set. If $\beta_k=1$ for some $k\le k^*$, then the interval $J_k$ is invariant, $\hat e_k$ is stable but not asymptotically stable fixed point and all points from
		$J_k\setminus \{e_k\}$ are eventually stable but not asymptotically stable fixed points or 2-periodic orbits.
		\end{remark}

\begin{theorem}\label{t3}
{\em Dynamics in each of the invariant segments $J_*$ and $J_k$ is defined by Theorems \ref{T2}
and \ref{Jan1}. Further, any  trajectory enters the union $\cup_{k\le k^*} J_k \cup J_*$
of these segments after a finite number of iterations.}
\end{theorem}

%\begin{remark}	\label{rJ}
%	kakie est attractors, kakie est unstable invariant sets -- varianty
%\end{remark}

\section{Proofs}

The proof proceeds in several steps and will be divided into a few lemmas.

Let us set %the set $\Omega$ as following
$$
\Omega=\{x\in J^*: \hbox{there is an } n\in{\mathbb N} \hbox{ such that } T^n(x)=0\}, \qquad \Sigma=\overline{\Omega}.
$$
By definition, the sets $\Omega$ and $\Sigma$ are invariant for $T$ and, moreover, if $y\in J^*$ and $T(y)\in\Omega$, then $y\in\Omega$.

\begin{lemma}\label{l1}
{\em   Suppose that a map $F$ is linear on segments $[d,c]$ and $[c,e]$, and the absolute value of the slope of $F$ on these segments is $\alpha$ and $\beta$, respectively. Suppose that
 \begin{equation}\label{ine2}
  \alpha^{-1}+\beta^{-1}<1.
  \end{equation}
   Then, for any $\lambda\in ( {\alpha}^{-1}+{\beta}^{-1},1)$ and any segment $[a,b]\subset [d,e]$, % the following inequality holds
  \begin{equation}\label{ine1}
    b-a<\lambda|F([a,b])|.
  \end{equation}
  }
\end{lemma}

%
%\begin{lemma}\label{l1}
%{\em   Suppose that
% \begin{equation}\label{ine2}
%  \alpha_k^{-1}+\beta_k^{-1}<1.
%  \end{equation}
%   Then, there is a $\lambda_k\in (0,1)$ such that for every segment $[a,b]\subset (q_{k+1},q_k)$, % the following inequality holds
%  \begin{equation}\label{ine1}
%    b-a<\lambda_k|T([a,b])|.
%  \end{equation}
%  }
%\end{lemma}

{\sf Proof.}
If  $b\leq c$, then $|F([a,b)])|=\alpha(b-a)$ and, similarly,
if $a\geq c$, then $|F([a,b)])|=\beta(b-a)$. Each of these equalities implies \eqref{ine1}.
%the proof is obvious.
Now, assume that $a<c<b$. Then,
$$
  |F([a,b])|\ge \max\{|F([a,c])|,|F([c,b])|\}=\max\{\alpha(c-a),\beta(b-c)\}.
$$
Without loss of generality we can assume that $\alpha(c-a)\geq\beta(b-c)$, hence
$$
  c-a=(b-a)-(b-c)\geq (b-a)-\frac{\alpha}{\beta}(c-a).
$$
%That is,
%$$
 % c-a\geq \frac{\beta}{\alpha+\beta}(b-a)
%$$
%and therefore,
Therefore,
$$
  |F([a,b])|\ge\alpha(c-a)\geq \frac{\alpha\beta}{\alpha+\beta}(b-a).
$$
%Since ${\alpha_k\beta_k}/{(\alpha_k+\beta_k)}>1$ for $\alpha_k^{-1}+\beta_k^{-1}<1$,
This implies  \eqref{ine1} for any $\lambda\in ( {\alpha}^{-1}+{\beta}^{-1},1)$. \hfill $\Box$

%\begin{remark}\label{remslope1}

%Since $\alpha_k$ and $\beta_k$ are increasing with $k$, we can assume that if \eqref{ine2} holds for $\alpha_{k_1}$ and $k_2$ and $k_1<k_2$, then $\lambda_{k_1}>\lambda_{k_2}$.
%{\bf D: It seems to be an extra assumption that $\alpha_k$ and $\beta_k$ are increasing?}
%\end{remark}

%There are 2 cases.
%
%\begin{itemize}
%
%\item
%$T(r_{k_0+1})>e_{k_0}$.
%In this case, the invariant segment $J$ attracts all the trajectories.
%
%\item
%$T(r_{k_0+1})<e_{k_0}$. In this case, the attractor has a more complicated structure.
%For this case, we need to establish a few facts, which are summarized in the following Lemmas \ref{l1}--\ref{l4}.
%
%\end{itemize}

\begin{lemma}\label{lJ1}
{\em $J^*$ is an invariant segment for the map $T$.}
\end{lemma}
{\sf Proof.}
First note that $[q_{k^*+1},r_{k^*}]\subset J^*$ and $T([q_{k^*+1},r_{k^*}])=J^*$. Since local maximum values $p_k$ of the map $T$ are decreasing with $k$, one has $\max_{x\in[0,q_{k^*}]}T(x)=p_{k^*}$. Hence, if $p_{k^*}\leq q_{k^*}$, then $T(J^*)=J^*$. On the other hand, if $q_{k^*}<p_{k^*}$, then
%since the graph of $T$ between $q_{k^*}$ and $T(r_{k^*})$ lies under or on the line $y=x$ then
the definition of $k^*$ implies that $T(x)\leq x \leq p_{k^*}$ for all $x\in [q_{k^*},p_{k^*}]$ (this part of the graph lies under the line $y=x$) and therefore $T(J^*)=J^*$ again.
 \hfill $\Box$

\begin{lemma}\label{lJ2}
{\em Suppose that $[0,\gamma]\subset J^*$. Then $[0,\gamma]\subset T([0,\gamma])$.}
\end{lemma}
{\sf Proof.}
Denote by $k_1$ the minimal number such that $e_{k_1}\leq \gamma$.
% and $\gamma<e_{k_1-1}$ if $e_{k_1-1}$ exists.
If $e_{k_1}=\gamma$ then the relation $T(r_{{k_1}+1})>e_{k_1}$ implies $[0,\gamma]\subset T([r_{{k_1}+1},q_{{k_1}+1}]) = T([0,\gamma])$.
 %Note that if $e_{k_1}=\gamma$, then $T(r_{{k_1}+1})\neq e_{k_1}$ because $[0,\gamma] \neq J^*$.
  If $e_{k_1}<\gamma\leq r_{k_1}$, then 
  %the point $(\gamma,T(\gamma))$ lies over the line $y=x$ 
  $\gamma<T(\gamma)$
  and therefore $[0,\gamma]\subset T([q_{{k_1}+1},\gamma])\subseteq T([0,\gamma])$. On the other hand, if $\gamma>r_{k_1}$, then there are two options. If $k_1>k^*$, then $\gamma <e_{k_1-1} < T(r_{k_1})$. If $k_1=k^*$, then $\gamma<T(r_{k_1})$ because $[0,\gamma] \neq J^*$. Hence, in both cases, $T([0,\gamma])=T([q_{k_1+1},r_{k_1}])=[0,T(r_{k_1})]\supset [0,\gamma]$.
 \hfill $\Box$

\begin{lemma}\label{l3}
{\em For every open interval $(a,b)$ satisfying $\Sigma\cap (a,b)\ne \emptyset$, there is a segment $[c,d]\subset(a,b)$ and an $n$ such that $T^{n}([c,d])=J^*$, $T^n(c)=0$, and $T^i$ is linear on the segment $[c,d]$ for each $i\leq n$.}
\end{lemma}

%{\bf D: V formulirovke lemmy $T(c)=0$ ili $T^k(c)=0$?}
%{\bf N: Done (Razumeetsya $T^n(c)=0$)}

{\sf Proof. }
Since $\Sigma=\overline{\Omega}$, there is a point $y\in \Omega$, $y\neq 0$ such that $y\in(a,b)$. Hence, there is an $n_1\in \mathbb{N}$ such that $T^{n_1}(y)=0$ and $T^{i}(y)\neq 0$ for $i<n_1$. The only condensation point of local extrema of the map $T$ is zero. Hence, if a point $z$ is a condensation point of local extrema of the iterated map $T^{r}$, then $T^{r-1}(z)=0$. Since $T^{i}(y)\neq 0$ for $i<n_1$, we conclude that
%without loss of generality we can assume that
there exists a point $l_1$ satisfying $a<y<l_1<b$ such that the map $T^{i}$ is linear on the segment $[y,l_1]$ for each $i\le n_1$.
%Since $T^{k_1}(y)=0$ and $T^{i}(y)\neq 0$ for $i<k_1$, it follows that there exists $[y,l_1]\subset[y,l]$ such that $T^{k_1}$ is linear on $[y,l_1]$.
  Denote by $\eta_1$ the number such that $T^{n_1+1}([y,l_1])=[0,\eta_1]$.

Since $T^{n_1}(y)=0$ and $T^{n_1}$ is linear on $[y,l_1]$, we have
\[
T^{n_1+1}(x)=T\left((x-y)\frac{T^{n_1}(l_1)}{l_1-y}\right)
\]
on $(y,l_1]$.
%then the graph of $T^{n_1+1}$ is piecewise linear on $[y+\epsilon,l_1]$ for any $0<\epsilon<l_1-y$. To be more precise, to obtain the graph of $T^{n_1+1}$ on $[y,l_1]$ one should restrict the graph of $T$ (consider $T$ from $0$ to $T^{n_1}(l_1)$) and then contract it.
Therefore, the graph of the iterated map $T^{n_1+1}$ is piecewise linear on $[y+\epsilon,l_1]$ for any $0<\epsilon<l_1-y$.
Hence, there exists a segment $[c_1,d_1]\subset[y,l_1]$ such that $T^{n_1+1}(c_1)=0$, $T^{n_1+1}(d_1)=\eta_1$ and $T^{n_1+1}$ is linear on $[c_1,d_1]$.

Similarly, we denote by $\eta_2$ the number such that $T^{n_1+2}([c_1,d_1])=T([0,\eta_1])=[0,\eta_2]$ and find the segment $[c_2,d_2]\subseteq[c_1,d_1]$ such that $T^{n_1+2}(c_2)=0$, $T^{n_1+2}(d_2)=\eta_2$, and $T^{n_1+2}$ is linear on $[c_2,d_2]$.
From Lemmas \ref{lJ1} and \ref{lJ2} it follows that $[0,\eta_1]\subseteq [0,\eta_2]$.

  %  Indeed, denote by $k_1$ a number such that $e_{k_1}\leq\eta_1$ but $e_{k_1-1}\nleq\eta_1$. If $e_{k_1}\leq\eta_1\leq r_{k_1}$ then $(\eta_1,T(\eta_1))$ lies over or on the line $y=x$ and then $[0,\eta_1]\subseteq T([q_{{k_1}+1},\eta_1])\subseteq T([0,\eta_1])=[0,\eta_2]$. In the other hand, if $\eta_1>r_{k_1}$ but $e_{k_1-1}\nleq\eta_1$ lets consider two cases. If $k_1<k^*$ then $\eta_1 <e_{k_1-1} \leq T(r_{k_1})$. Then obviously $T([0,\eta_1])=T([q_{k_1+1},r_{k_1}])=[0,T(r_{k_1})]$. If $k_1=k^*$ then $J^*=[0,T(r_{k_1})]$. Hence in both cases $[0,\eta_1]\subseteq T([0,\eta_1])=[0,\eta_2]$.

Continuing this line of argument, we obtain segments
$$
[c_1,d_1]\supseteq [c_2,d_2]\supseteq\cdots\supseteq [c_j,d_j]\supseteq\cdots
$$
and
$$
[0,\eta_1]\subseteq [0,\eta_2]\subseteq\cdots\subseteq [0,\eta_j]\subseteq\cdots
$$
such that   $T^{n_1+j}([c_j,d_j])=[0,\eta_{j}]$, $T^{n_1+j}$ is linear on $[c_j,d_j]$, $T^{n_1+j}(c_j)=0$ and $T^{n_1+j}(d_j)=\eta_j$. Note also that $T([0,\eta_j])=[0,\eta_{j+1}]$.

Let us consider the limit $\eta=\lim_{j\rightarrow\infty}\eta_j$. Note that we can fix a small $\epsilon_1>0$ such that $T([\epsilon_1,\eta_j])=[0,\eta_{j+1}]$ for all $j$. Then, by continuity, $T([0,\eta])=[0,\eta]$. From Lemma \ref{lJ2} it follows that there are no invariant subsegments containing $0$ inside $J^*$.
Hence $[0,\eta]=J^*$. Fix an $n$ such that $\eta_{n}>r_{k^*}$. Taking into account that $J^*=[0,p_{k^*}]$, %{\bf (``and $T$ decreases on $[r_{k^*},p_{k^*}]$" -- D: we don't need this)}
we see that $T([0,\eta_{n}])=J^*$.
Therefore, the segment $[c_{n+1},d_{n+1}]$ satisfies $T^{n_1+n+1}([c_{n+1},d_{n+1}])=J^*$, $T^{n_1+n+1}(c_{n+1})=0$, and $T^i$ are linear mappings on $[c_{n+1},d_{n+1}]$ for $i\leq n_1+n+1$. \hfill $\Box$

\begin{lemma}\label{l4}
{\em Suppose that $T^{n}(\hat{x})=\hat{x}\in J^*$ and there is a segment $[a,\hat{x}]$ such that $T^{n}(a)=0$ and $T^{n}$ is linear on $[a,\hat{x}]$. Then $\hat{x}\in\Sigma$.}
\end{lemma}

{\sf Proof.}
Since $T^{n}(a)=0$, $T^{n}(\hat{x})=\hat{x}$, $\hat{x}>a$ and $T^{n}$ is linear on $[a,\hat{x}]$, there is a unique $a_1\in[a,\hat{x}]$ such that $T^n(a_1)=a$. Obviously, $T^{2n}(a_1)=0$, $T^{2n}(\hat{x})=\hat{x}$, and $T^{2n}$ is linear on $[a_1,\hat{x}]$. Arguing in a similar way, we obtain a sequence $a_1<a_2<a_3<\cdots$ such that
$T^{(i+1)n}(a_i)=0$ (hence, $a_i\in\Omega$) and $T^{(i+1)n}$ is linear on $[a_i,\hat{x}]$. Since the slope of the graph of $T^{(i+1)n}$ on $[a_i,\hat{x}]$ tends to infinity, it follows that $a_i\rightarrow \hat{x}$. Hence $\hat{x}\in\Sigma$.
\hfill $\Box$

\begin{remark}\label{remeksigma}
From Lemma \ref{l4} it follows that $e_k\in \Sigma$ for $k\geq k^*$.
\end{remark}

\begin{lemma}\label{l5}
{\em Let $\alpha^{-1}_{k^*}+\beta^{-1}_{k^*}>1$. Suppose that $x\in J^*$, $x\notin \Sigma$, but $T(x)\in \Sigma$. Then, $x$ is a local maximum point of the map $T$ and there exists an $n\in\mathbb{N}$ such that
$T^{n}(x)=e_{k^*}$.}
\end{lemma}

{\sf Proof. }
From Remark \ref{remslope2} it follows that $J_{k^*}$ is invariant. By definition of $\Sigma$ and Remark \ref{remeksigma}, this implies $J_{k^*}\cap\Sigma=e_{k^*}$.

Since $x\notin \Sigma$ and $\Sigma$ is closed, there is a neighborhood $U$ of the point $x$ such that $U\cap\Sigma=\emptyset$. Since there are no points from $\Omega$ inside $U$, it follows that $T^{i}$ is continuous on $U$ and $\Omega\cap T^{i}(U)=\emptyset$  for every $i$. The relationships $T(x)\in\Sigma=\overline{\Omega}$ and $\Omega\cap T(U)=\emptyset$ imply that $T(x) \notin {\rm Int}\,(T(U))$. Furthermore, taking into account that $x\neq 0$
and neither the right end point $T(r_{k^*})$ of $J^*$ nor its iterations belong to $\Sigma$,
%{\bf (D: pochemu net?)} %{\bf (N: Potomu chto ona lezhit v invariantnom otrezke, v kotorom tol'ko odna tochka, $e_{k^*}$, lezhit v $\Sigma$, no v etu tochku nikto iz etogo otrezka ne perehodit)},
we conclude that either $T(y)\leq T(x)$ for all $y\in U$ or $T(y)\geq T(x)$ for all $y\in U$. Since all local minimum points of $T$ inside $J^*$ belong to $\Omega$, we see that $x$ is a local maximum point for $T$, and there exists a $b>0$ such that
the interval $\theta:=[T(x)-b,T(x)]$ satisfies $\theta\cap\Omega=\emptyset$.

Since $f_k<e_{k-1}<T(r_k)$ for $k\geq k^*+1$ (cf.~\eqref{fk}), Remark \ref{remslope2} implies that $\alpha^{-1}_k+\beta_k^{-1}<1$ for $k\geq k^*+1$. Hence, from Lemma \ref{l1} and the assumption that $\alpha_k$ and $\beta_k$ are increasing %and Remark \ref{remslope1}
 it follows that there is a $\lambda<1$ such that if $T^{\bar{n}}(\theta)\subset[0,e_{k^*}]$, $T^{\bar{n}}(\theta)\cap\Omega=\emptyset$, then
\begin{equation}\label{n1}
  |T^{\bar {n}}(\theta)|<\lambda|T^{\bar{n}+1}(\theta)|.
\end{equation}
This inequality %\eqref{n1}
implies that there is an ${n}'$ such that if $T^i(\theta)\subset[0,e_{k^*}]$ for $i\leq{n}'$, then
$T^{{n}'}(\theta)\cap\Omega\neq\emptyset$, which contradicts  $\theta\cap\Omega=\emptyset$. Consequently, there is an $\tilde{n}$ such that $T^{\tilde{n}}(\theta)\cap {\rm Int}(J_{k^*})\neq\emptyset$. On the other hand,
$e_{k^*}\in \Sigma$ implies $e_{k^*} \notin {\rm Int}\, (T^{\tilde{n}}(\theta))$ because $\theta\cap \Omega=\emptyset$. %Hence, $T^{\tilde n}(\theta)\cap\Omega=\emptyset$. 
Therefore, the segment $T^{\tilde{n}}(\theta)$ satisfies
 $T^{\tilde{n}}(\theta)\subseteq J_{k^*}$. Finally, since $\Sigma$ is invariant and $T(x)\in \Sigma$ by assumption, $T^{\tilde{n}}(x)\in \Sigma\cap J_{k^*}=e_{k^*}$. \hfill $\Box$

% {\bf D: Sleduyushii Remark vrode tavtologiya: esli funktsiya strogo monotonna na intervale, to tam net lokalnih maximumov. Ty chto hotel skazat?}
%
%\begin{remark}\label{remmaximum}
%It is easy to see that if $f$ is linear on $(a,b)$ then there are no local maximum points of $f$ inside $(a,b)$.
%\end{remark}

%{\bf D: Pravilno li zamenit v sleduyushem Ramark $i\le n$ na $0\le i\le n-1$?}

\begin{remark}\label{remsigma}
Lemma \ref{l5} implies that if $\alpha^{-1}_{k^*}+\beta^{-1}_{k^*}>1$, $x \notin \Sigma$, and $T^{i}(x)$ is not a local maximum point of $T$ for $i\leq n-1$, then $T^{n}(x)\notin \Sigma$.
\end{remark}

\begin{remark}\label{remab1}
A slight modification of the proof of Lemma \ref{l5} shows that if $\alpha^{-1}_{k^*}+\beta^{-1}_{k^*}=1$, $x\in J^*$ and $x\notin \Sigma$ but $T(x)\in \Sigma$, then either $x$ is a local maximum point of $T$ or $x=T(r_{k^*})$. In both cases, there exists an $n\in\mathbb{N}$ such that
$T^{n}(x)=e_{k^*}$.
\end{remark}

\begin{remark}\label{remab2}
  From the proof of the Lemma \ref{l5} it follows that if $\alpha^{-1}_{k^*}+\beta^{-1}_{k^*}\geq1$, then for every point $x\in J^*\backslash\Sigma$ there
  exists an $n$ such that $T^i(x)\in J_{k^*}$ for $i\geq n$. Indeed, since $\Sigma$ is closed, there is a segment $\vartheta\subset J^*$ such that $x\in\vartheta$ and $\vartheta\cap\Omega=\emptyset$. Therefore, arguing in the same way as in the proof of Lemma \ref{l5}, we obtain that there is an $n$ such that $T^{n}(\vartheta)\subseteq J_{k^*}$. Since $J_{k^*}$ is invariant for $\alpha^{-1}_{k^*}+\beta^{-1}_{k^*}\geq1$, one has $T^{i}(\vartheta)\subseteq J_{k^*}$ for all $i\geq n$.
\end{remark}

\begin{lemma}\label{l6}
{\em Suppose that $\alpha_{k^*}^{-1}+\beta_{k^*}^{-1}<1$. Then, $\Sigma=J^*$.}
\end{lemma}

{\sf Proof. }
Suppose that there is an $x\in J^*$ such that $x \notin \Sigma$. Since $\Sigma$ is closed, there is a segment $\theta\subset J^*$ such that $x\in\theta$ and $\theta\cap\Omega=\emptyset$.
%Since $\alpha_{k}^{-1}+\beta_k^{-1}<1$ for $k\geq k^*$ then from
From Lemma \ref{l1} and the assumption that $\alpha_k$ and $\beta_k$ are increasing it follows that there is a $\mu<1$ such that if %$T^{\bar{n}}(\theta)\subset J^*$,
$T^{{n}}(\theta)\cap\Omega=\emptyset$, then
\begin{equation}\label{n2}
  |T^{{n}}(\theta)|<\mu|T^{{n}+1}(\theta)|.
\end{equation}
As the segment $J^*$ is invariant,
from \eqref{n2} it follows that there exists an $\tilde{n}$ such that $T^{\tilde{n}}(\theta)\cap\Omega\neq\emptyset$. But this contradicts the fact that $\theta\cap\Omega=\emptyset$.
Hence, we conclude that every $x\in J^*$ belongs to $\Sigma$.
\hfill $\Box$

\begin{theorem}\label{T1}
{\em $\Sigma$ is a chaotic invariant set.}
\end{theorem}

{\sf Proof. }
There are two cases, when $J_{k^*}$ is invariant and $J_{k^*}$ is not invariant. We present a proof for the more complicated
case when $J_{k^*}$ is invariant. The other case, when $J_{k^*}$ is not invariant, can be treated similarly.

First, let us prove sensitive dependence on initial conditions and density of periodic points in $\Sigma$.
 Denote $\zeta={e_{k^*}}/{3}$. Consider any point $x\in\Sigma$ and its neighborhood $(a,b)\ni x$. By Lemma \ref{l3}, there is a segment $[c,d]\subset(a,b)$ and an $n\in \mathbb{N}$ such that $T^{n}([c,d])=J^*$, $T^n(c)=0$, and $T^i$ is linear on the segment $[c,d]$ for each $i\leq n$. Hence, there is a point $\tilde{x}\in[c,d]$ such that $T^n(\tilde{x})=\tilde{x}$, and from Lemma \ref{l4} it follows that $\tilde{x}\in\Sigma$. Also,
 $T^{n}([c,d])=J^*$ implies that there is a point $z\in(c,d)$ such that $T^n(z)=e_{k^*}$. Since $T^i$ is linear on the interval $(c,d)$ for each $i\leq n$, it follows from
 %Remarks \ref{remmaximum} and
 Remark \ref{remsigma} that $z\in\Sigma$.
% {\bf D: nado kak-to bolee podrobno ob'yasnit, chto "$T^i$ is linear on the segment $[c,d]$ for each $i\leq n$" $\Rightarrow$ "$T^i (z)$ is not a local maximum point of $T$ for $i \le n$"} {\bf N: Ya dobavil remark.}
 Obviously,
\begin{equation}\label{zeta}
\max\{|T^n(x)-T^n(z)|,|T^n(x)-T^n(c)|\}\ge |T^n(z)-T^n(c)|/2=e_{k^*}/2>\zeta.
\end{equation}
Since the interval $(a,b)\ni x$ can be chosen arbitrarily small,
the relationships $\tilde x, c, z\in (a,b)\cap \Sigma$, $T^n(\tilde x)=\tilde x$ and \eqref{zeta} prove
sensitivity to initial conditions and density of periodic points in $\Sigma$.

It remains to prove the existence of a dense orbit in $\Sigma$.

Note that since $J_{k^*}$ is invariant, $\Omega\subset [0,e_{k^*}]$.

For any $n\in \mathbb{N}$, let us consider a collection of segments $\{I^i_n\}$, $i=1,...,l(n)$, with $l(n)\leq n$, $I^i_n\subseteq J^*$ such that %{\bf (D: Some problem with $I_n^i$ here)}
\begin{itemize}

\item $|I^i_n|\leq {e_{k^*}}/{n}$, $i=1,...,l(n)$;

\item $I^i_n\cap\Omega\neq\emptyset$, $i=1,...,l(n)$;

\item $\Sigma\subset \bigcup_{i=1}^{l(n)} I^i_n$.

\end{itemize}
Let us consider all the segments $I^i_n$, $i=1,...,l(n)$, $n\in \mathbb{N}$, and number them as follows:
$$
L_1=I_1^1;
$$
$$
L_{r+1}=I^{i+1}_n \quad \hbox{ if }\quad L_r=I^i_n, \ i<l(n);
$$
$$
L_{r+1}=I^1_{n+1} \quad \hbox{ if }\quad L_r=I^{l(n)}_n.
$$

From Lemma \ref{l3} it follows that there is a segment $E_1\subset L_1$ and an $n_1$ such that $T^{n_1}(E_1)=L_2$ and $T^i$ is linear on the segment $E_1$ for each $i\leq n_1$. Denote $H_1=E_1$. Similarly, there is a segment $E_2\subset L_2$ and an $n_2$ such that $T^{n_2}(E_2)=L_3$ and $T^i$ is linear on the segment $E_2$ for each $i\leq n_2$. Hence, there is a segment $H_2\subset H_1$ such that $T^{n_1+n_2}(H_2)=L_3$ and $T^i$ is linear on the segment $H_2$ for each $i\leq n_1+n_2$. Continuing in a similar fashion
, we obtain a sequence of nested segments
$$
H_1\supseteq H_2\supseteq...\supseteq H_r\supseteq...
$$
and a sequence $n_1, n_2,...,n_r,...$ such that
$$
T^{n_1+n_2+...+n_r}(H_r)=L_{r+1}
$$
and $T^j$ is linear on the segment $H_r$ for each $j\leq  n_1+n_2+...+n_r$.

Denote $D_r=H_r\cap\Sigma$. As an intersection of two closed sets, $D_r$ is closed.
Consider the non-empty intersection
$$
D=\bigcap_{r=1}^{\infty}D_r
$$
of the closed nested sets $D_r$. Take any point $x\in D$. By construction, the forward orbit of $x$ is dense in $\Sigma$.
\hfill $\Box$

Combining Remark \ref{remab2}, Lemma \ref{l6}, and Theorem \ref{T1}, we obtain Theorem \ref{T2}.

 \begin{remark}\label{Jan2} Arguing in a same fashion as in the proof of Theorem \ref{T1} we obtain the following statement.
		 	Suppose that $F:[c,d]\to[c,d]$ is continuous and for any interval $[a,b]\subseteq[c,d]$ there is an $n>0$
		 	such that $F^n([a,b])=[c,d]$. Then, $[c,d]$ is a chaotic invariant set for $F$.
		 \end{remark}

\medskip
\noindent
{\sf Proof of Theorem \ref{Jan1}. }

Cases $\beta_k<1$ and $\beta_k=1$ for $k\leq k^*$ are trivial.

\begin{figure}[H]
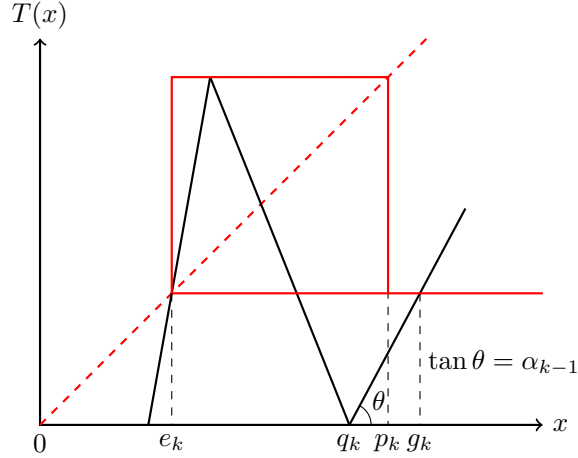

 	 \centering
 	 \includestandalone[width=0.5\textwidth]{techn}
 	 \caption{Restriction of the map $T$ to an interval $J_k=[e_k,p_k]$, cf.~Fig.~\ref{fig:types}(c).\label{fig:techn}}
 	 \end{figure}

\medskip
Now consider the case $\alpha_{k}^{-1}+\beta_{k}^{-1}<1$ for some $1\leq k\leq k^*-1$. Denote
\[
g_k=\min\{x> f_k:\ T(x)=e_k\},\qquad 1 \leq k\le k^*,
\] 
with $f_k$ defined by \eqref{fk}.
  Since $\alpha_{k-1}={e_k}/({g_k-q_k})$ (see Fig.~\ref{fig:techn}), we have ${e_k}/{\alpha_{k-1}}=g_k-q_k$. Hence from \eqref{star} it follows that $g_k>p_k$.
 The conclusion of the theorem in the case $g_k>p_k$
is well known and follows from the general theory of unimodal maps (see for example \cite{devaney}).

%Hence, the proof for the case $\alpha_{k}^{-1}+\beta_{k}^{-1}<1$  is complete.

\medskip
The last case  $\alpha_{k}^{-1}+\beta_{k}^{-1}>1$ with $\alpha_k, \beta_k>1$ will be considered by induction.

\medskip

Let us introduce the sequences defined by the recurrent relations
  \[
  \xi_0=\alpha_k, \ \ \nu_0=\beta_k;\qquad \xi_{i+1}=\nu_i^2, \ \ \nu_{i+1}=\xi_i\nu_i,\qquad i=0,1,2,\ldots
  \]
 By assumption, $\xi_0^{-1}+\nu^{-1}_0>1$. Further, since $\alpha_k,\beta_k>1$, the sequence
 $\xi_i^{-1}+\nu_i^{-1}$ monotonically decreases to zero. Therefore, there is an index $j=j(\alpha_k,\beta_k)\ge 1$ such that
 \[
 \begin{array}{ccc}
 \xi_i^{-1}+\nu_i^{-1}>1 & {\rm for} & 0\le i\le j-1,\\
  \xi_i^{-1}+\nu_i^{-1}\le 1 & {\rm for} & i\geq j.
 \end{array}
  \]
 We will show that the statement of the theorem is true for $J_k$ with $N=j-1$ using the induction in $j$.

\medskip

{\em Basis. } We start from the case $j(\alpha_k,\beta_k)=1$. In other words,
$\beta_{k}^{-1}(\alpha_{k}^{-1}+\beta_{k}^{-1})\le 1$. 
For the sake of brevity, the following argument is conducted under the assumption that the strict inequality 
\begin{equation}\label{aa}
\beta_{k}^{-1}(\alpha_{k}^{-1}+\beta_{k}^{-1})< 1
\end{equation}
%We will refer to the following well known fact, which can be proved using the same line of argument as in the proof of Theorem \ref{T1} above.
holds. The case of the equality can be done by a slight modification of the same argument.

 From $\beta_k>1$ and $\alpha_{k}^{-1}+\beta_{k}^{-1}> 1$ it follows that the interval $J_k$
 and its subinterval $G_k$ defined by \eqref{JG} are invariant for $T$, and each trajectory starting from $J_k\setminus\{e_k\}$ enters the interval $G_k$ after a finite number of iterations.
 %, see Fig.~\ref{??}. {\bf D: What should be on this figure?}

 Consider the map $T^2: J_k\to J_k$.
 This is a piecewise linear map with 4 linear segments, i.e.~there is a partition
 \[
 e_k<u_k<r_k<v_k<p_k
 \]
 of the segment $J_k$ such that
 \[
 \frac{d}{dx}T^2(x)=\left\{
 \begin{array}{rl}
 \alpha_k^2, & x\in(e_k,u_k),\\
 - \alpha_k\beta_k, & x\in(u_k,r_k),\\
  \beta_k^2, & x\in(r_k,v_k),\\
   - \alpha_k\beta_k, & x\in(v_k,p_k).
 \end{array}
 \right.
  \]
  Notice that the function $T^2$ reaches its maximum value $p_k$ at the points $x=u_k,v_k$ and has a local minimum at the point $x=r_k$.
 % which has the value $T(r_k)$ at its maximum point(s)  and the value $T^2(r_k)$ at its minimum point lying in ${\rm int}\, G_k=(T^2(r_k),T(r_k))$, see Fig.~\ref{??}. If $T^2$ has 3 linear segments, then their slopes are
 %$-\alpha_k\beta_k$, $\beta^2_k$, $-\alpha_k\beta_k$. If $f$ has 4 linear segments, their slopes are $\alpha^2_k$, $-\alpha_k\beta_k$, $\beta^2_k$, $-\alpha_k\beta_k$.

   Consider an arbitrary segment $\Delta\subseteq G_k$.
Let us show that $\hat e_k\in {\rm Int}\, (T^\ell(\Delta))$ for some $\ell$.
Assume the contrary. Then, any iteration $T^i(\Delta)$ contains at most one of the
extremum points $u_k, r_k,v_k$.
% Now, suppose that  $\hat e_k\not\in {\rm int}\, T^\ell(\Delta)$ for all $\ell$.
 By Lemma \ref{l1}, relation \eqref{aa} implies that there is a $\lambda_1>1$
such that if $u_k\not \in\Delta$, then
 $|T^2(\Delta)|>\lambda_1 |\Delta|$. On the other hand,
 if $u_k \in\Delta$, then $T^2(\Delta)=[\delta,p_k]$ with $\delta>v_k$ because
 $\hat e_k\in {\rm Int}\,(T^4([v_k,p_k]))$. Hence, $T^4$ is piecewise linear on $\Delta$
 with the two slopes $-\alpha_k^3\beta_k$ and $\alpha_k^2\beta_k^2$.
 From \eqref{aa} and $\alpha_k>1$ it follows that
   \[
  \alpha_k^{-2}\beta_k^{-1}(\alpha_{k}^{-1}+\beta_{k}^{-1})<1,
  \]
hence by Lemma \ref{l1} there is a $\lambda_2>1$ such that
  $|T^4(\Delta)|>\lambda_2 |\Delta|$. Thus, either   $|T^2(\Delta)|>\lambda_1 |\Delta|$ or  $|T^4(\Delta)|>\lambda_2 |\Delta|$ (or both).  Continuing the iterations, we see that there exists $\lambda$ such that for any $n$ there should be
  an $n_1$  such that $|T^{n_1}(\Delta)|>\lambda^n|\Delta|$, which contradicts the
  invariance of $G_k$.

 We conclude that $\hat e_k\in {\rm Int}\, (T^\ell(\Delta))$ for some $\ell$, which immediately  implies that
   $T^m(\Delta)=G_k$ for some $m>\ell$.
   Now, the conclusion of the theorem with $N=0$, $A_1={\rm Int}(G_k)$ and $\Lambda_k=G_k$ follows from Remark \ref{Jan2}.

  \medskip

  {\em Induction step.}

  Assume that the conclusion of the theorem is valid 
  %for all $\alpha_k$, $\beta_k$ such that 
  for all $j\leq j_0$.
  Now, assume that $j(\alpha_k,\beta_k)=j_0+1$. Consider the segment $J_k'=[\hat e_k,p_k]\subset J_k$.
  Since $j(\alpha_k,\beta_k)=j_0+1\ge 2$, we have $\beta_{k}^{-1}(\alpha_{k}^{-1}+\beta_{k}^{-1})>1$, hence $J_k'$ is invariant for $T^2$.
  Further,
    \begin{equation}\label{bb}
    T(J_k')\cap J_k'=\{\hat e_k\}.
    \end{equation}
    Moreover,
  the restriction of $T^2$ to $J_k'$ has the same
  shape as the restriction of $T$ to $J_k$, i.e.~$T^2$ is piecewise linear on $J_k'$ with two slopes
  $\xi_1=\beta_k^2$ and $-\nu_1=-\alpha_k\beta_k$.
  Clearly, $j(\xi_1,\nu_1)=j(\alpha_k,\beta_k)-1$. Therefore,
  by the induction assumption, the conclusion of the theorem holds for the restriction of $T^2$ to $J_k'$.
  Combining this statement with relation \eqref{bb} and the fact that $\hat e_k$ is an unstable fixed point of $T$,
  we obtain the statement of the theorem for $j=j_0+1$ and formulas \eqref{sigma}. $\Box$

  \medskip
  \noindent
  {\sf Proof of Theorem \ref{t3}. }

  For each $x\in J$ such that $x \notin \cup_{k\le k^*} J_k \cup J_*$ define $1\leq k(x)\leq k^*$ such that $k(x)=\min\{k : e_k < x\}$. Since the graph of $T$
 can  lie over the line $y=x$ only at points which belong to the union $\cup_{k\le k^*} J_k \cup J_*$, then, after finitely many iterations, $x$ is either mapped to $\cup_{k\le k^*} J_k \cup J_*$ or to a point $y$ such that $k^*\ge k(y)>k(x)$. Thus, after finitely many iterations, $x$ is mapped to $\cup_{k\le k^*} J_k \cup J_*$. $\Box$

 %It follows from \eqref{stent} that
 %\[
 %\frac1\beta=\frac1b-1,\qquad \frac1\alpha=\frac{2b-1}{1-a},
 %\]
 %hence
 %\[
 %\frac1{\beta^2}+\frac1{\alpha\beta}=1+\frac{2b-1}{b^2}\left(\frac{1-b}{1-a}-1\right)
 %\]
 %and since $0<b<a<1$ and $b>1/2$, we see that
 %\[
 %\frac1{\beta^2}+\frac1{\alpha\beta}<1.
 %\]
 %Further, $\alpha>\beta$ implies that
 %\[
 %\frac1{\alpha^2}+\frac1{\alpha\beta}<1.
 %\]

 \section{Two-dimensional system with PWL saturation function}
 In this section, we apply the results of Section \ref{main} to the system
 \begin{equation}
 \label{2D}
 	\begin{cases}
 	x_{n+1}=\lambda x_n+ (\sigma-\lambda)s_n,\\
 	s_{n+1}=\Phi(s_n+x_{n+1}-x_n)
 	\end{cases}
 	\end{equation}
 	with $n\in \mathbb{N}_0$, where $\Phi$ is the PWL {\em saturation} function 
 	\begin{equation}\label{saturation}
 	\Phi(x)=\left\{\begin{array}{rlc}
 	-1 & {\rm if} & x\le -1,\\
 	x & {\rm if} & |x|<1,\\
 	1 & {\rm if} & x\ge 1.
 	\end{array}
 	\right.
 	\end{equation}
 	The phase space for this system is the horizontal strip
 	\begin{equation}\label{LLLL}
 	Q=\left\{(x,s): \ x \in \mathbb{R}, \ -1\le s \le 1 \right\}.
 	\end{equation}
 	We will use the short notation
 	\[
 	(x_{n+1},s_{n+1})=F(x_n,s_n)
 	\]
 	for \eqref{2D}. By definition, the function $F$ maps $Q$ into itself.
 	
 	Below, we consider the domain of parameters
 	\[
 	\Pi=\{(\sigma,\lambda): \ \sigma>1,\  -1<\lambda<0\}.
 	\]
 	Outside this domain, the global attractor of system \eqref{2D} consists either of equilibrium points or a period 2 orbit, or a union thereof \cite{siamds}. In $\Pi$, dynamics are more interesting.
 	
 	%For a large parameter domain, it has been shown that the global attractor of this map
 	%can consist of a segment of, or just two, equilibrium points, a period two orbit bifurcating from
 	%the segment of equilibrium points via a degenerate flip border collision bifurcation, or a union of equilibrium points and period two orbit(s). Further, in the complementary domain of parameters, the map has been shown to possess infinitely many periodic orbits of which at most one is stable.
 	%The proof of this statement, pointing to the likely presence of chaos, has been based on the reduction to a Poincar\'e section.

\begin{figure}[H]
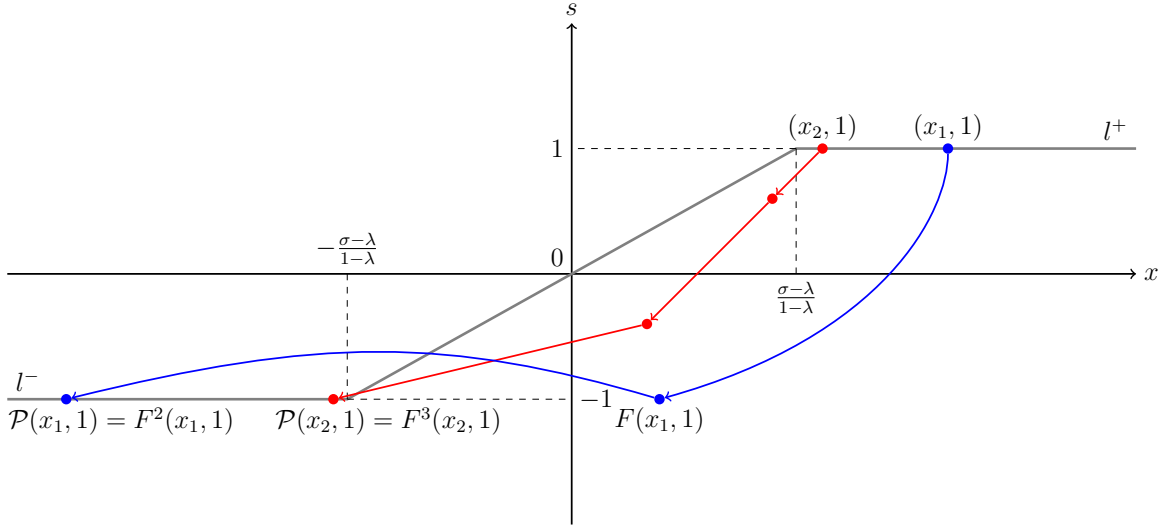

 	\begin{center}
 	\includestandalone[width = \textwidth]{poincare}
 	\end{center}
 	\caption{First-hitting map $\mathcal P: l^+ \to \overline{l^-}$ defined by trajectories of the map $F$. The trajectory starting from $(x_1, 1)\in l^+$ arrives at the half line $\overline{l^-}$ after two iterations under the map $F$. The trajectory starting from $(x_2, 1)\in l^+$ arrives at $\overline{l^-}$ after three iterations.  \label{fig:poincare}
 			}
 	\end{figure}

 	It is easy to see that equilibrium points of system \eqref{2D}
 	form the segment
 	\begin{equation}\label{eq3}
 	%EF
 	E^-E^+=\left\{(x,s): \  x=\frac{(\sigma-\lambda) s}{1-\lambda},\ -1\leq s\leq 1\right\}\subset Q.
 	\end{equation}
 	Further, consider the two horizontal half-lines starting from
  	the end points $E^+$ and $E^-$ of this segment,
 	\[
 	l^+=\left\{(x,s): x> \frac{\sigma-\lambda}{1-\lambda}, \ s=1\right\},
 	\qquad
 	l^-=\left\{(x,s): x< -\frac{\sigma-\lambda}{1-\lambda}, \ s=-1\right\}.
 	\]
 	It has been shown in \cite{siamds} that
 	any trajectory of \eqref{2D} starting from the half-line $l^+$ arrives at the closed half-line $\overline{l^-}$ after  finitely many
 	iterations; note that similarly trajectories starting on $l^-$
 	reach $\overline{l^+}$ because the map $F$ is odd.
 	Hence, we can define the first-hitting map ${ \mathcal P}: l^+\to
 	\overline{l^-}$
 	as ${\mathcal P}(x,s)=F^k(x,s)$ where $F^{k}(x,s)\in \overline{l^-}$ and
 	$F^{i}(x,s)\not
 	\in \overline{l^-}$ for $i=1,\ldots,k-1$, see Fig.~\ref{fig:poincare}.
 	This map can be represented by the scalar function
 	$f: (\frac {\sigma-\lambda}{1-\lambda},\infty)\to [\frac {\sigma-\lambda}{1-\lambda},\infty)$ defined by the
 	formula
 	\begin{equation}\label{TP}
 	(-f(x),-1)={\mathcal P}(x,1),\qquad x\in\left(\frac
 	{\sigma-\lambda}{1-\lambda},\infty\right).
 	\end{equation}
 	
 	%Using the fact that the map $f$ is odd and the rays $l_1\cup E$ and $l_3\cup F$ are centrally symmetric, we identify these rays by the negative identity map. With this notation, and setting by definition ${T}(E)=E$, ${T }$ will be considered as the map of $l_1 \cup E$ into itself.
 Further, the function $f$ is piecewise linear on every interval $[a,b]\subset  (\frac {\sigma-\lambda}{1-\lambda},\infty)$ with local minimum points
 	\begin{equation}\label{qk}
 	\hat q_k=\frac{-\frac{1}{1-\lambda}-\frac{1-\sigma^k}{1-\sigma}}{\frac{\sigma^k}{\sigma-\lambda}-\frac{1-\sigma^k}{1-\sigma}},\qquad f(\hat q_k)=\frac {\sigma-\lambda}{1-\lambda},\quad k=1,2,\ldots,
 	\end{equation}
 and local maximum points
 \begin{equation}\label{rk}
 \hat r_k=\frac{2+(\sigma-\lambda)\frac{1-\sigma^k}{1-\sigma}}{1-\sigma^k+(\sigma-\lambda)\frac{1-\sigma^k}{1-\sigma}},\qquad \hat{p}_k  = f(\hat r_k)=\sigma-\lambda(\hat{r}_k-1),\quad k=1,2,\ldots,
 \end{equation}
 where
 \[
 \hat q_1> \hat r_1>\hat q_2>\hat r_2>\cdots;\qquad \hat q_k, \hat r_k\to \frac{\sigma-\lambda}{1-\lambda};
 \]
also, $f'(x)=-\lambda$ for $x>\hat q_1$ \cite {siamds}. These relations imply that
\[
\lim_{k\to \infty}f(\hat r_k) = %\frac{ -2\lambda(\sigma-1)}{1-\lambda}
\sigma+\lambda -\lambda \frac{\sigma-\lambda}{1-\lambda}>\frac{\sigma-\lambda}{1-\lambda}=\lim _{k\to\infty} \hat r_k
\]
(recall that $(\sigma,\lambda)\in \Pi$), hence $f$ has an essential discontinuity at the left end of its domain $ (\frac {\sigma-\lambda}{1-\lambda},\infty)$, and therefore there is a $k^{**}$ defined by
\[
k^{**}=\min \{k\ge 1: f(\hat r_k) \ge \hat r_k \}.
\]

 	It is convenient to shift the origin and consider
 	the function $T: \mathbb{R}_+\to \mathbb{R}_+$ given by
 	\begin{equation}\label{2DT}
 	T(0)=0;\qquad
 	T(x)=f\left(\frac{\sigma-\lambda}{1-\lambda}+x\right)-\frac{\sigma-\lambda}{1-\lambda}, \quad x>0.
 	\end{equation}
The above properties of $f$ imply that $T$ is a ``saw map'' with the sequences \eqref{pqr} defined by
 	\begin{equation}\label{newseq}
 %	q_k=\frac{-\frac{1}{1-\lambda}-\frac{1-\sigma^k}{1-\sigma}}{\frac{\sigma^k}{\sigma-\lambda}-\frac{1-\sigma^k}{1-\sigma}}-\frac{\sigma-\lambda}{1-\lambda},\qquad r_k=\frac{2+(\sigma-\lambda)\frac{1-\sigma^k}{1-\sigma}}{1-\sigma^k+(\sigma-\lambda)\frac{1-\sigma^k}{1-\sigma}}-\frac{\sigma-\lambda}{1-\lambda}
 q_k=\hat q_{k+k^{**}-1}-\frac{\sigma-\lambda}{1-\lambda},\qquad
  r_k=\hat r_{k+k^{**}-1}-\frac{\sigma-\lambda}{1-\lambda},\qquad
   p_k=\hat p_{k+k^{**}-1}-\frac{\sigma-\lambda}{1-\lambda}
 \end{equation} for $k\ge 1$. If $k^{**} > 1$ then $r_0$ and $p_0$ can be also defined by \eqref{newseq}. If $k^{**}=1$ then we can put $r_0$ to be any number greater than $q_1$ such that $p_0 = f(r_0) - \frac{\sigma-\lambda}{1-\lambda}>p_1$.
 	
 	%(and $T(x)<x$ for $x>q_1$).
 %	and
% 	\[
 %	p_k=-\lambda\left(r_k+\frac{2(\sigma-1)}{1-\lambda}\right).
 %	\]
 	Further, the slopes of $T$ on the segments $[q_{k+1},r_k]$ and $[r_k,q_k]$, respectively, are given by
 	\[
 	 \alpha_k=-\left(\sigma^{k+1}-(\sigma-\lambda)\frac{1-\sigma^{k+1}}{1-\sigma}\right), \qquad
 	 -\beta_k=-\lambda\left(\sigma^k-(\sigma-\lambda)\frac{1-\sigma^k}{1-\sigma}\right).
 	\]
 	%and also $T'(x)=-\lambda$ for $x>q_1$.
 	It is easy to see that both sequences $\alpha_k$ and $\beta_k$ are increasing
 	and $\alpha_k>\beta_k$ for all $k$.
 	%Since $p_k\to \frac{ -2\lambda(\sigma-1)}{1-\lambda}>0$ as $k\to \infty$, the map $T$ has an essential discontinuity at zero.
 	Also note that according to \eqref{rk}, all the points $(r_k,p_k)$ lie on a straight line $l$ with the slope $-\lambda$.
 	Moreover, the map \eqref{2DT} satisfies technical condition \eqref{star} as a consequence of the fact that $\alpha_k > \alpha_{k-1}$ and the following proposition.
 %	{\bf D: We need to renumber $p_k,q_k,r_k$ so that the first maximum satisfying $p_k>r_k$ becomes $p_1$. We also assume in Section 3 that $q_1>p_1$ --- is it true here?}

 \begin{proposition}\label{p1}
 	The map \eqref{2DT} satisfies
 	\begin{equation*}\label{techn}
 	q_k+\frac{e_{k}}{\alpha_{k}}>T(r_k)
 	\end{equation*}
 	for all $1\le k<k^*$.
 	\end{proposition}
 	
\begin{figure}[h]
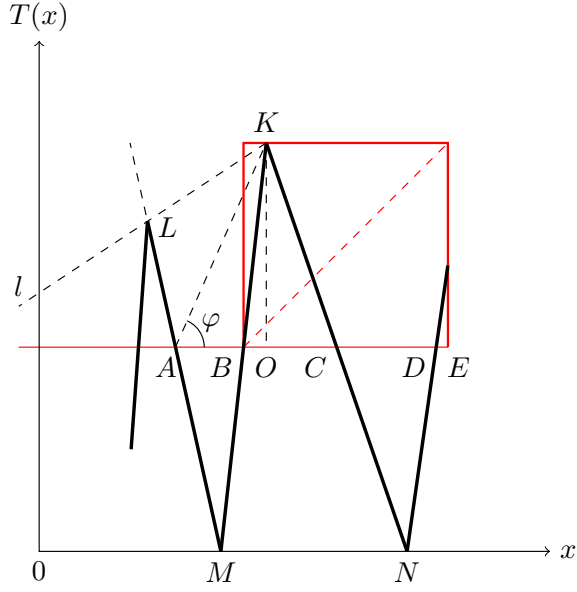

 	 \centering
 	 \includestandalone[width=0.5\textwidth]{pic}
 	 \caption{Proof of Proposition \ref{p1}. The thick polyline is the graph of $T$. Coordinates of the points are as follows: $N = (q_k, 0)$, $M = (q_{k+1}, 0)$, $K = (r_k, p_k)$, $L = (r_{k+1}, p_{k+1})$, $B=(e_k, e_k)$, $E = (p_k, e_k)$. The angle $\varphi=\angle{KAB}$ satisfies $\tan{\varphi} > 1$.\label{figABCDE}}
 	 \end{figure}

 	 \noindent
 	 {\sf Proof.} 
 	 Suppose that the inequality is not valid for some $k$. Since the sequences $\alpha_k$ and $\beta_k$ increase, Fig.~\ref{figABCDE} implies the following inequlities:
 	 \[
 	 |AO|<|AC|<|MN|<|BD|\le |BE|=|OK|.
 	 \]
 	 Hence, the straight line
 	 the straight line $AK$ has a slope greater than 1. Therefore, the straight line
 	 $l$ with the slope $-\lambda<1$ passing through the point $K$ intersects the line $MA$ above the horizontal line $AB$. That is, the intersection point $L=(r_{k+1},p_{k+1})$ satisfies
 	 $p_{k+1}>e_k$. By definition of $k^*$, this implies $k\ge k^*$. \hfill $\Box$
 	
 	\medskip
 	For a given pair of parameters $(\sigma,\lambda)\in\Pi$, the above explicit formulas for $q_k,p_k,r_k,\alpha_k,\beta_k$ allow us to find the intervals $J^*$ and $J_k$
 	for the map $T=T_{\sigma,\lambda}$
 	and to determine which case of Theorems \ref{T2} and \ref{Jan1} applies to each of these intervals
 	depending on whether $\alpha_k^{-1}+\beta_k^{-1}$ is less or greater than 1 and whether $\beta_k$
 	is less or greater than 1. We performed these computations numerically in the rectangle $1<\sigma<3$, $0>\lambda>-1$
 	%with the step $\delta\lambda=...$, $\Delta\sigma=...$
 	at $1000\times 1000$ points. 
 	
	For the sake of brevity we will use the following classification of the intervals $J_k$. If $\beta_k \le 1$, then we say that the interval $J_k$ is of type I; if $\beta_k > 1$ and $\alpha_k^{-1} + \beta_k^{-1} \ge 1$, then $J_k$ is of type II; and, if $\beta_k > 1$ and $\alpha_k^{-1} + \beta_k^{-1} < 1$, then $J_k$ is of type III (see Fig.~\ref{fig:types}(a-c)).
	
\begin{figure}[h]
	\begin{center}
	\includegraphics[width=0.5\textwidth]{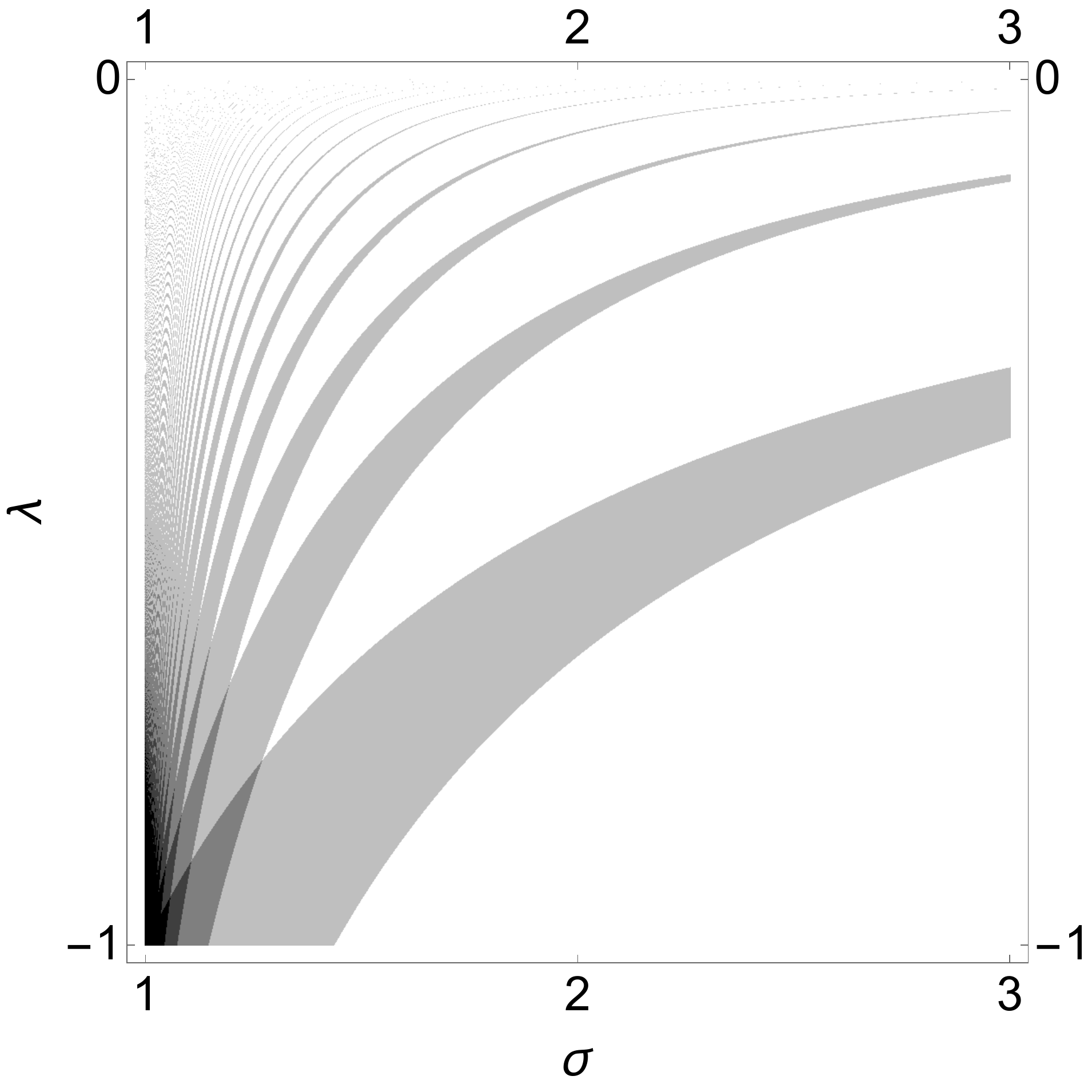}
	\caption{The number $k^*(\sigma,\lambda)$ of intervals $J_k$ is represented according to the following convention: white corresponds to $k^*=1$, light gray to $k^*=2$, gray to $k^*=3$, dark gray to $k^*=4$, black to $k^*\ge5$. The number $k^* $ increases as $\lambda \to -1$ and $\sigma \to 1$.} \label{fig2}
	\end{center}
	\end{figure}

	Our numerical findings can be summarized as follows. First, the number of intervals $J_k$ increases  as $\lambda\to -1$, $\sigma\to 1$. In other words, $k^* \to \infty$ as $\lambda\to -1$, $\sigma\to 1$. Fig.~\ref{fig2} shows the value of $k^*=k^*(\sigma,\lambda)$.
		
\begin{figure}[H]
	\begin{center}
	\subfloat[\label{fig3a}]{\includegraphics[width=0.49\textwidth]{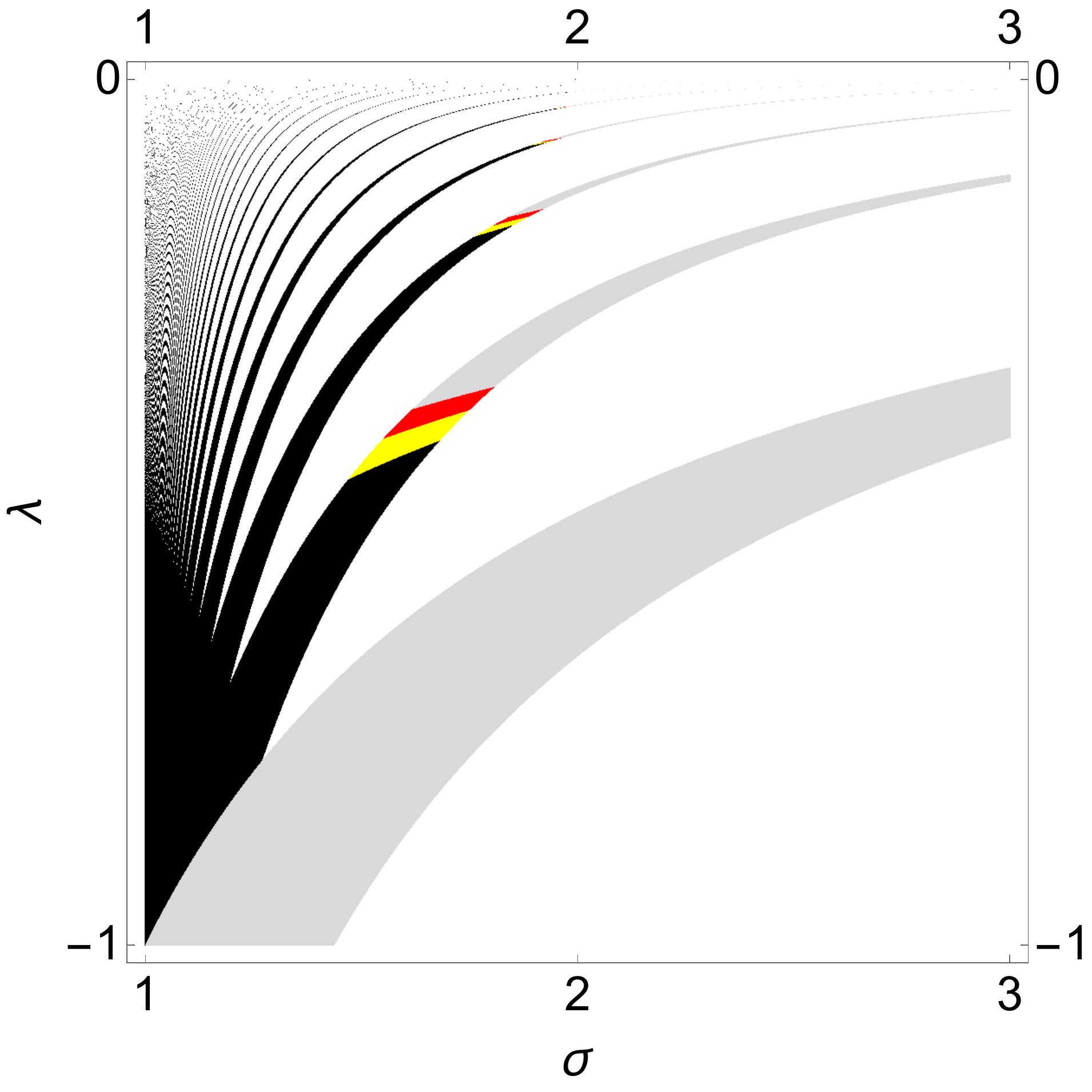}}
%	\hspace{40pt}
	\subfloat[\label{fig3b}]{\includegraphics[width=0.49\textwidth]{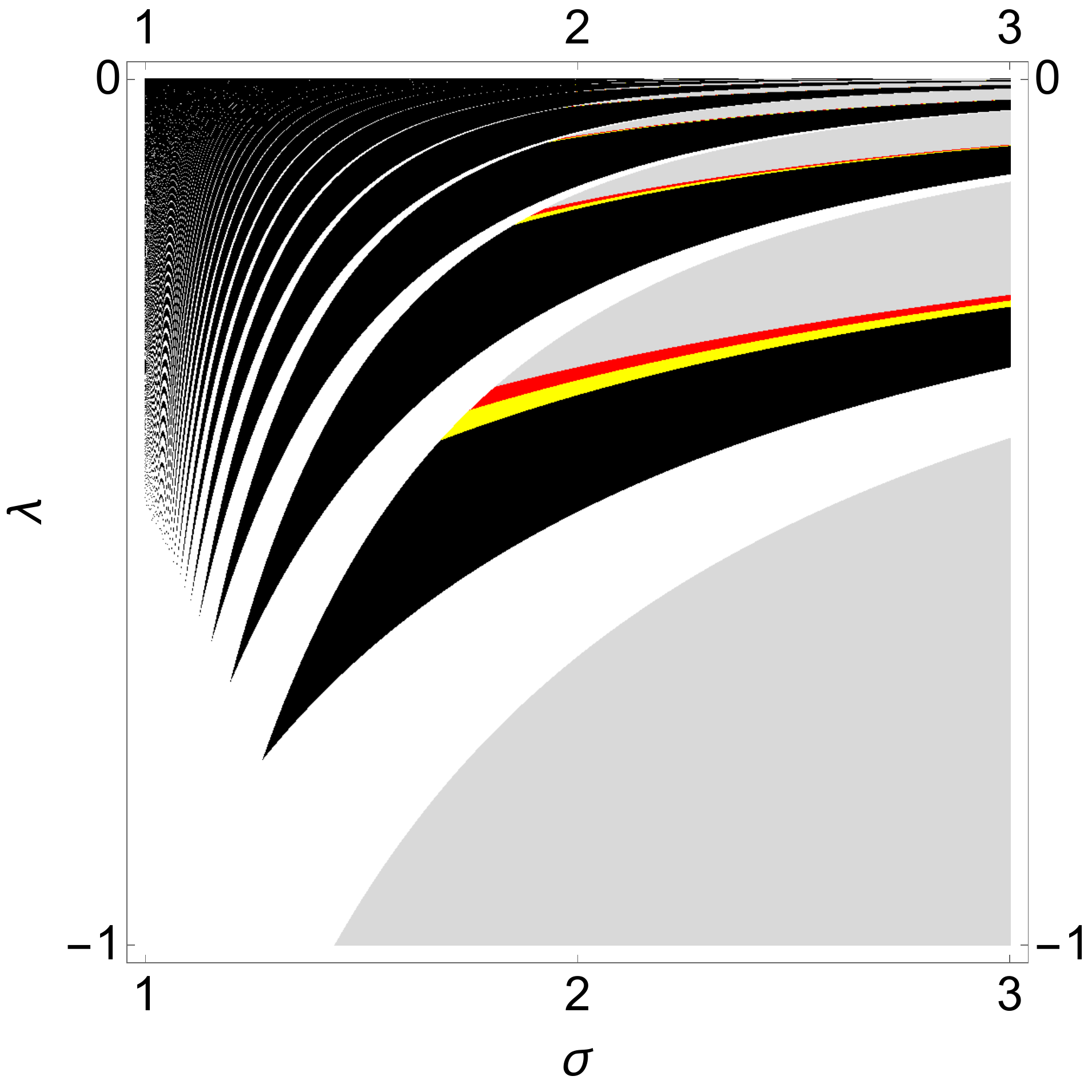}}
	\caption{The type of the interval $J_1$ depending on the values of the parameters $\lambda$ and $\sigma$.  $(a)$ $k^* \ge 2$ ($J_1 \not\subset J^*$); $(b)$: $k^* = 1$ ($J_1 \subset J^*$). We adopt the following convention: light gray corresponds to type I, yellow to type II with the chaotic attractor $\Lambda_k\subset J_k$ consisting of one segment, red to type II with the chaotic attractor $\Lambda_k\subset J_k$ consisting two disjoint segments, black to type III, white on panel $(a)$ to $k^*=1$, white on panel $(b)$ to $k^*\geq 2$.} \label{fig3}
	\end{center}
\end{figure}

	Second, only the rightmost segment $J_1$ can have any of the types I, II or III. All the other segments (if they exist) are of type III. This observation suggests that the map \eqref{2DT} does not  have segments of types I and II simultaneously, {\it i.e.} the asymptotically stable fixed point of the map \eqref{2DT} does not coexist with a chaotic invariant set of type \eqref{Lambda} composed of a finite number of closed intervals. Fig.~\ref{fig3}(a) shows the type of the interval $J_1$ (depending on the values of parameters $\lambda$ and $\sigma$) when $k^* \ge 2$. In this case,  $J_1 \not\subset J^*$. Similarly, Fig.~\ref{fig3}(b) shows the type of the segment $J_1$ for the case $k^*=1$ when $J_1 \subset J^*$.
	
	Third, we found that if the rightmost interval $J_1$ is of type II, then the chaotic attractor $\Lambda_k$ contained in this interval is either a segment or a union of two disjoint segments (cf.~\eqref{Lambda}).

	Note that the domain ${\mathcal D}$ of existence for the stable fixed point of $T$ (the rightmost segment $J_1$ is of type I) has been found in \cite{siamds}:
\[
{\mathcal D}=\bigcup_{i=1}^\infty \left\{(\lambda,\sigma):\
\dfrac{\sigma^i-1}{\sigma-1}\le-\frac1{\lambda}<\sigma^i, \ \sigma>1, \
-1<\lambda<0\right\}.
\]
	
	\begin{remark}
		Any point $(x,1)\in l^+$ with $x\in (\hat q_{k+1}, \hat q_k)$ reaches the half-line $l^-$ in $k+1$ iterations under the map \eqref{2D} \cite{siamds}. This implies that every periodic orbit of $T$ 
		%that belongs to $J_k$ 
		corresponds to a periodic orbit of the map \eqref{2D}. In particular, since $J_k\subset (q_{k+1},q_k)$, the fixed points $e_k$, $\hat e_k$ correspond to $(2(k+1))$-periodic orbits of \eqref{2D}. Stable orbits of $T$ correspond to stable orbits of \eqref{2D}.
		
%Comment that there are $k$ iterations (?) to reach from $l^+$ to $l^-$ if $x\in (q_{k+1}, q_k)$.
%Interpret the above bullet points in terms of dynamics using Theorems 1 and 2.
\end{remark}

\begin{remark}
	System \eqref{2D} has been alternatively interpreted in \cite{siamds} as follows. 	
Let $s_0\in[-1,1]$ and let $\{x_n\}$, $n\in \mathbb{N}_0$, be a real-valued sequence.
Consider
the sequence
$$
s_{n+1}=\Phi(s_n+x_{n+1}-x_n),\quad n\in\mathbb{N}_0,
$$
with the saturation function \eqref{saturation}. The mapping of the pair 
 $s_0,\, \{x_n\}$ to the sequence $\{s_n\}$
according to this formula is known in different disciplines\footnote{The stop operator with discrete time inputs/outputs is used, for example, in applications to economics \cite{Gocke,networks}. Applications in engineering and physics typically use the continuous time extension of this operator.} as the {\em stop} operator $\mathcal{S}$ \cite{krasnoselskii,brokate},
Prandtl's model of elastic-ideal plastic element \cite{prandtl}, the backlash nonlinearity, and the one-dimensional Moreau sweeping process \cite{Kunze1,moreau}.
 Here $s_0$ is called the initial state, $\{x_n\}$ is called the input, $
 \{s_n\}={\mathcal S}[s_0,\{x_n\}]
 $ is called the output (or, the variable state) of the stop operator.
 %, and it is convenient to represent ${\mathcal S}$ by the input-output diagram shown in Fig.~\ref{fig:0d}.
  With this notation, system \eqref{2D} is equivalent to
 the simple feedback loop coupling a linear unit with the stop operator:
 \[
 \{x_{n+1}\}=\lambda \{x_n\} +(\sigma-\lambda){\mathcal S}[s_0,\{x_n\}],
 \]
 which is the interpretation used in \cite{siamds}.
 \end{remark}

\section*{Acknowledgments}
The authors thank P. Gurevich for a stimulating discussion of the results. D. R.
and P. K. acknowledge the support of NSF through grant DMS-1413223. N. B. acknowledges Saint-Petersburg State University (research grant 6.38.223.2014), the Russian Foundation for Basic Research (project No 16 -- 01 -- 00452) and DFG project SFB 910. The publication was financially supported by the Ministry of Education and Science of the Russian Federation (the Agreement number 02.A03.21.0008).

\bibliographystyle{abbrv}
\bibliography{biblio}

\end{document}

%% file: saw_pic.tex
\begin{tikzpicture}[scale=1.2]
\tikzset{
  font={\fontsize{12pt}{12}\selectfont}}
	\draw [<->,thick] (0,8) node (yaxis) [above] {$T(x)$}
		|- (12,0) node (xaxis) [right] {$x$};
	\coordinate (o) at (0, 0);
	\draw (o) node[anchor=north] {0};
	\draw[dashed, red, thick, name path = diag] (0, 0) -- (8, 8);
	\coordinate (q1) at (10.5, 0);
	\coordinate (p1) at (5.5, 6.5);
	\coordinate (q2) at (3.3, 0);
	\coordinate (p2) at (2.2, 3.8);
	\coordinate (q3) at (1.4,0);
	\coordinate (p3) at (1, 3);
	\coordinate (q4) at (0.7, 0);
	\coordinate (p5) at (0.4, 2.5);
	\coordinate (q6) at (0.3, 0);
	\coordinate (p6) at (0.21, 2);
	\coordinate (q7) at (0.17, 0);
	\coordinate (p7) at (0.14, 1.8);
	\draw[thick] (q1) node [anchor=north] {$q_k$}--(p1); 
	\draw[thick, name path = p1q2] (p1) -- (q2);
	\draw[thick] (q2) -- (p2);
	\draw[thick, name path = p2q3] (p2)-- (q3);
	\draw[thick] (q3)--(p3);
	\draw[thick] (p3) -- (q4);
	\draw[thick] (q4) -- (p5) -- (q6) -- (p6) -- (q7)--(p7);
	\path  [name intersections={of=diag and p1q2, by=e1}];
	\path  [name intersections={of=diag and p2q3, by=e2}];
	\draw[dashed] (p1)--($(q1)!(p1)!(q2)$) node [anchor = north] {$r_k$};
	\draw[dashed] (p2)--($(q2)!(p2)!(q3)$) node [anchor = north] {$r_{k^*}$};
	\node at (q2) [anchor = north] {$q_{k^*}$};
	\draw[dashed] let \p{p2}=(p2) in (p2) -- (\y{p2}, \y{p2}) -- (\y{p2}, 0) node [anchor=north] {$p_{k^*}$};
	\draw[thick, red] let \p{e2}=(e2), \p{p2}=(p2) in (e2) -- (\x{e2}, \y{p2})--(\y{p2}, \y{p2})--(\y{p2}, \y{e2})--(e2); 
	\draw[decoration={brace, mirror, raise=32pt, amplitude=15pt}, decorate, thick] let \p{p2}=(p2) in (0,0) -- (\y{p2}, 0) node[pos=0.5,anchor=north,yshift=-46pt] {$J^*$};
	\draw[dashed] let \p{e2}=(e2) in (e2) -- (\x{e2}, 0) node [anchor=north] {$e_{k^*}$};
	\draw[decoration={brace, mirror, raise=15pt, amplitude=8pt}, decorate, thick, red] let \p{e2}=(e2), \p{p2}=(p2) in (\x{e2}, 0)--(\y{p2}, 0) node[pos=0.5,anchor=north,yshift=-22pt] {$J_{k^*}$};
	\draw[dashed] let \p{e1}=(e1) in (e1) -- (\x{e1}, 0) node [anchor=north] {$e_{k}$};
	\draw[dashed] let \p{p1}=(p1) in (p1) -- (\y{p1}, \y{p1}) -- (\y{p1}, 0) node [anchor=north] {$p_{k}$};
	\draw[thick, red] let \p{e1}=(e1), \p{p1}=(p1) in (e1) -- (\x{e1}, \y{p1})--(\y{p1}, \y{p1})--(\y{p1}, \y{e1})--(e1);
	\draw[decoration={brace, mirror, raise=15pt, amplitude=8pt}, decorate, thick, red] let \p{e1}=(e1), \p{p1}=(p1) in (\x{e1}, 0)--(\y{p1}, 0) node[pos=0.5,anchor=north,yshift=-22pt] {$J_{k}$};
	\draw[thick] (q1)--(12, 2.7);
\end{tikzpicture}

%% file: type1.tex
\begin{tikzpicture}[scale=5]
\coordinate (o) at (0, 0);
\coordinate (T) at (1/3,1);
\draw[red, thick] (o) node [anchor=north, black] {$(e_k, e_k)$} -- (0, 1) -- (1, 1) -- (1, 0) -- (o);
\draw[red, dashed, thick, name path = diag] (o) -- (1, 1);
\draw[black, thick] (0, 0) -- (T) node [anchor=south] {$\left(r_k,p_k\right)$};
\draw[black, thick, name path = TE] (T) -- (1, 17/25);
\path  [name intersections={of=diag and TE, by=Ehat}];
\node at (Ehat) [anchor=east] {$\left(\hat{e}_k, \hat{e}_k\right)$};
\end{tikzpicture}

%% file: type2.tex
\begin{tikzpicture}[scale=5]
\coordinate (o) at (0, 0);
\coordinate (T) at (1/3,1);
\draw[red, thick] (o) node [anchor=north, black] {$(e_k, e_k)$} -- (0, 1) -- (1, 1) -- (1, 0) -- (o);
\draw[red, dashed, thick, name path = diag] (o) -- (1, 1);
\draw[black, thick] (0, 0) -- (T) node [anchor=south] {$\left(r_k, p_k\right)$};
\draw[black, thick,  name path = TE] (T) -- (1, 1/4);
\path  [name intersections={of=diag and TE, by=Ehat}];
\node at (Ehat) [anchor=east] {$\left(\hat{e}_k, \hat{e}_k\right)$};
\end{tikzpicture}

%% file: type3.tex
\begin{tikzpicture}[scale=5]
\coordinate (o) at (0, 0);
\coordinate (T) at (1/3,1);
\draw[red, thick] (o) node [anchor=north, black] {$(e_k, e_k)$} -- (0, 1) -- (1, 1) -- (1, 0) -- (o);
\draw[red, dashed, thick, name path = diag] (o) -- (1, 1);
\draw[black, thick] (0, 0) -- (T) node [anchor=south] {$\left(r_k, p_k\right)$};
\draw[black, thick,  name path = TE] (T) -- (4/5 ,0);
\path  [name intersections={of=diag and TE, by=Ehat}];
\node at (Ehat) [anchor=east] {$\left(\hat{e}_k, \hat{e}_k\right)$};
\end{tikzpicture}

%% file: techn.tex
\begin{tikzpicture}[scale=1]
	\draw [<->,thick] (0,5) node (yaxis) [above] {$T(x)$}
		|- (6.5,0) node (xaxis) [right] {$x$};
	\coordinate (o) at (0, 0);
	\draw (o) node[anchor=north] {0};
	\draw[dashed, red, thick, name path = diag] (0, 0) -- (5, 5);
	\coordinate (p1) at (5.5, 2.8);
	\coordinate (qk) at (4, 0);
	\coordinate (pk) at (2.2, 4.5);
	\coordinate (q3) at (1.4,0);
	\draw[thick, name path = p1qk] (p1) -- (qk);
	\draw[thick] (qk) -- (pk);
	\draw[thick, name path = pkq3] (pk)-- (q3);
	\path  [name intersections={of=diag and pkq3, by=ek}];
	\node at (qk) [anchor = north] {$q_{k}$};
	\draw[dashed] let \p{pk}=(pk) in (pk) -- (\y{pk}, \y{pk}) -- (\y{pk}, 0) node [anchor=north] {$p_{k}$};
	\draw[thick, red] let \p{ek}=(ek), \p{pk}=(pk) in (ek) -- (\x{ek}, \y{pk})--(\y{pk}, \y{pk}) -- (\y{pk}, \y{ek});
%	\draw[thick, red, name path = sq] let \p{ek}=(ek), \p{pk}=(pk) in (\y{pk}, \y{ek})--(ek);
	\draw[thick, red, name path = sq] let \p{ek}=(ek), \p{pk}=(pk) in (6.5, \y{ek})--(ek);
	\path  [name intersections={of=sq and p1qk, by=gk}];
	\draw[dashed] (gk)--($(o)!(gk)!(xaxis)$) node (gk0) [anchor=north] {$g_k$};
	\begin{scope}
	\path[clip] (qk) -- (gk) -- (gk0);
	\draw[black] (qk) circle (8pt);
	\end{scope}
	\node at ($(qk)+(40:5mm)$) {$\theta$};
	\node at (6, .8) {$\tan{\theta} = \alpha_{k-1}$};
	\draw[dashed] (ek) -- ($(o)!(ek)!(xaxis)$) node [anchor=north] {$e_k$};
\end{tikzpicture}

%% file: poincare.tex
\tikzset{
  font={\fontsize{12pt}{12}\selectfont}}
\begin{tikzpicture}[scale = 2]
\pgfmathsetmacro{\A}{1.78947}
	\coordinate (o) at (0, 0);
	\node  at (o)  [anchor = south east]{$0$};
	\draw[thick, ->] (-4.5, 0) -- (4.5, 0) node (xaxis) [anchor=west] {$x$};
	\draw[thick, ->] (0, -2) -- (0, 2) node (saxis) [anchor=south] {$s$};
	\coordinate (l1) at (\A, 1);
	\coordinate (l2) at (-\A, -1);
	\draw[very thick, gray] (-4.5, -1) node [anchor=south west, black] {$l^-$} -- (l2) -- (l1) -- (4.5, 1) node [anchor=south east, black] {$l^+$};
	\draw[dashed] (l1) -- ($(o)!(l1)!(saxis)$) node [anchor=east] {$1$};
	\draw[dashed] (l2) -- ($(o)!(l2)!(saxis)$) node [anchor=west] {$-1$};
	\draw[dashed] (l1) -- ($(o)!(l1)!(xaxis)$) node [anchor=north] {$\frac{\sigma-\lambda}{1-\lambda}$};
	\draw[dashed] (l2) -- ($(o)!(l2)!(xaxis)$) node [anchor=south] {$-\frac{\sigma-\lambda}{1-\lambda}$};
	\node[red, circle, fill=red, scale = 0.5] (a1) at (2, 1) {};
	\node[red, circle,fill=red, scale = 0.5] (a2) at (1.6, 0.6){};
	\node[red, circle, fill=red, scale = 0.5] (a3) at (0.6, -0.4){};
	\node[red, circle, fill=red, scale = 0.5] (a4) at (-1.9, -1){};

	\node[blue, circle, fill=blue, scale = 0.5] (b1) at (3, 1) {};
	\node[blue, circle,fill=blue, scale = 0.5] (b2) at (0.7, -1){};
	\node[blue, circle, fill=blue, scale = 0.5] (b3) at (-4.03, -1){};
	
	\draw[red, thick, ->] (a1) node [anchor = south, black] {$(x_2, 1)$}--(a2);
	\draw[red, thick, ->] (a2)--(a3);
	\draw[red, thick, ->] (a3)--(a4) node [anchor = north, black, xshift = 25pt] {$\mathcal{P}(x_2, 1)=F^3(x_2, 1)$};

	\draw[blue, thick, ->] (b1) node [anchor=south, black] {$(x_1, 1)$}.. controls (3, .5) and (2.5, -0.5) ..(b2) node [anchor=north, black] {$F(x_1, 1)$};
	\draw[blue, thick, ->] (b2) .. controls (-1, -0.5) and (-2, -0.5) .. (b3) node [black, anchor = north, xshift = 25pt] {$\mathcal{P}(x_1, 1)=F^2(x_1, 1)$} ;
\end{tikzpicture}

%% file: pic.tex
\begin{tikzpicture}[scale=2.5]
\coordinate (B) at (1, 1);
\coordinate (K) at (10/9,2);
\coordinate (N) at (9/5, 0);
\coordinate (M) at (8/9 ,0);
\coordinate (E) at (2, 1);
\coordinate (P) at (4/9, 2);
\draw[red, thick] (B) node [anchor=north east, color=black] {$B$} -- (1, 2) -- (2, 2) -- (E) node[anchor = north west, color=black, xshift=-4pt] {$E$};
\draw[red, name path = BE] (E) -- (-0.1, 1);
\draw[red, dashed] (B) -- (2, 2);
\draw[black, very thick] (M) node [anchor=north] {$M$} -- (K);
\draw[black, very thick, name path=KN] (K) node [anchor=south] {$K$} -- (N) node [anchor=north] {$N$};
\draw[black, very thick, name path = ND] (N) -- (2, 7/5);
\draw[black, dashed, name path=MP] (M) -- (P);
%\draw[black] (0, 0) -- (2, 0);
\node (o) at (0, 0) [anchor=north] {$0$};
\draw [<->, black] (-0,2.5) node (yaxis) [above] {$T(x)$} |- (2.5,0) node (xaxis) [right] {$x$};
\draw[black, dashed] (K) -- ($(B)!(K)!(E)$) node[anchor=north] {$O$};
\path [name intersections={of=BE and ND, by = D}];
\node at (D) [anchor = north east] {$D$}; 
\path [name intersections={of=BE and KN, by = C}];
\node at (C) [anchor = north east] {$C$};
\path [name intersections={of=BE and MP, by = A}];
\node at (A) [anchor = north east, xshift=4pt] {$A$};
\draw[black, dashed, name path=KL] (K) -- (-0.1, 6/5) node [above] {$l$};
\path [name intersections={of=KL and MP, by=L}];
\node at (L) [right, yshift=-2] {$L$};
\draw[black, very thick] (M)--(L);
\draw[black, dashed] (A) -- (K);
\draw[black, very thick] (L) -- (.45,0.5);
\begin{scope}
\path[clip] (A) -- (K) -- (B);
\draw[black] (A) circle (4pt);
\end{scope}
\node at ($(A)+(35:6pt)$) {$\varphi$};
\end{tikzpicture}